





\magnification\magstep1
\baselineskip =18pt
\overfullrule =0pt
\pageno=0

\def\reel{{\rm I}\!{\rm R}}
\def\comp{\;{}^{{}_\vert}\!\!\!{\rm C}}
\def\nat{{{\rm I}\!{\rm N}}}

\def\exp{{{\rm I}\!{\rm E}}}
\def\free{{\rm I}\!{\rm F}}
\def\1{1\mkern -4mu{\rm l}}
\def\n{\noindent}
\def\qed{\qquad {\rm qed}}
\def\pf{ \medskip \n {\bf Proof.~~}}

\def\M{{\cal M}}
\def\H{{\cal H}}
\def\C{{\cal C}}
\def\F{{\cal F}}
\def\S{{\cal S}}
\def\E{{\cal E}}
\def\A{{\cal A}}
\def\N{{\cal N}}

\def\BMO{{{\cal B}{\cal M}{\cal O}}}

\def\vp{\varepsilon}

\def\O{\Omega}
\def\o{\omega}
\def\f{\varphi}
\def\qed{{\hfill{\vrule height7pt width7pt depth0pt}\par\bigskip}}
\def\ie{{\it i.e.\/}\ }
\def\cf{{\it cf.\/}\ }

\centerline{\bf Non-commutative martingale inequalities}
\bigskip

\setbox1=\vtop{\hsize=6cm
\centerline {Gilles Pisier~*}
\centerline{Texas A\&M University and}
\centerline{Universit\'e Paris VI}
\centerline{Equipe d'Analyse; Case 186}
\centerline{4, Place Jussieu}
\centerline{75252 Paris Cedex 05; France}
\centerline{ gip@ccr.jussieu.fr}}

\setbox2=\vtop{\hsize=6cm
\centerline{Quanhua Xu}
\centerline{Universit\'e de Franche-Comt\'e}
\centerline{Equipe de Math\'ematiques}
\centerline{UFR des Sciences et Techniques}
\centerline{16, Route de Gray}
\centerline{25030 Besan\c con Cedex; France}
\centerline{ qx@math.univ-fcomte.fr}}
$$\line{\box1\hfill\box2}$$\footnote{~}{*Partially supported by the NSF}
\bigskip

\noindent{\bf Abstract.} We prove the analogue of the classical
Burkholder-Gundy inequalites for
non-commutative  martingales.   As applications we give a characterization
for an Ito-Clifford
integral to be an $L^p$-martingale via its integrand, and then extend the
Ito-Clifford integral theory
in $L^2$, developed by Barnett, Streater and Wilde, to $L^p$ for all
$1<p<\infty$. We include an
appendix on the non-commutative analogue of the classical Fefferman duality
between
$H^1$ and $BMO$.
\bigskip
\noindent{\bf Plan.}

\noindent 0. Introduction

\noindent 1. Preliminaries

\noindent 2. The main result

\noindent 3. Examples

\noindent 4. Applications to the Ito-Clifford integral

\noindent~~ Appendix \vfill\eject

\centerline{\bf  0. Introduction}

Recently, non-commutative (=quantum) probability theory has developed
considerably. In particular,
all sorts of non-commutative analogues of Brownian motion and martingales
have been studied following
the basic work of Parthasarathy and Schmidt. We refer the reader to P. A.
Meyer's exposition ([M])
and to the proceedings of the successive conferences on quantum probability
[AvW] for more details
and references. There are also intimate  connections with Harmonic Analysis
(cf. e.g. [Mi]).

Motivated by quantum physics, and after the pioneer works of Gross (cf.
[Gr1-2]), a Fermionic version
of Brownian motion and stochastic integrals was developed (see [BSW1]),
and  the optimal
hypercontractive inequalities have been finally proved ([CL]).

In this paper we will prove the non-commutative analogue of the classical
Burkholder-Gundy
inequalities from martingale theory.  We should point out that what
follows was originally inspired by some
recent work of Carlen and Kr\'ee, who had considered Fermionic versions of
the Burkholder-Gundy
inequalities. They obtained the inequality in Theorem 4.1 below in some
special cases, as well as
some sufficient conditions for the convergence of stochastic integrals in
the case $p\le 2$ (see
section 4 below for more on this).

One interesting feature of our work, is that the square function is defined
differently (and it must
be changed!) according to $p<2$ or $p>2$. This surprising phenomenon was
already discovered  by F.
Lust-Piquard in [LP] (see also [LPP])  while
 establishing  non-commutative versions of Khintchine's inequalities.

Let us briefly describe  our main inequality. Let $\M$ be a finite  von
Neumann algebra with a
normalized  normal faithful trace
$\tau$, and $(\M_n)_{\n\ge 0}$ be an increasing filtration of von Neumann
subalgebras of $\M$.
 Let $1<p<\infty$ and $(x_n)$ be a  martingale with respect to 
$(\M_n)_{\n\ge 0}$ in the usual
$L^p$-space $L^p(M,\tau)$ associated to $(M,\tau)$. Set
$d_0=x_0,\ d_n=x_n-x_{n-1}$. Then our main result reads as follows. If
$p\ge 2$, we have (with
equivalence constants depending only on $p$)
$$\sup_n\|x_n\|_p \approx \max\big\{
\| (\sum d_n^*d_n)^{1/2} \|_p,\  \| (\sum d_nd_n^*)^{1/2}
\|_p\big\}.\leqno(0.1)$$ This is no longer valid for $p<2$; however 
for $p<2$ the
``right" inequalities are
$$\sup_n\|x_n\|_p \approx
\inf\big\{\| (\sum a_n^*a_n)^{1/2} \|_p+ \| (\sum b_nb_n^*)^{1/2}
\|_p\big\},\leqno(0.2)$$ where the infimum
runs over all decompositions
$d_n=a_n+b_n$ of $d_n$ as a sum of martingale difference sequences 
adapted to the same filtration.

In particular, this applies to martingale transforms: given a martingale
 $(x_n)$ as above and an adapted bounded sequence $\xi=(\xi_n)$,  
\ie  such that 
$\xi_n\in \M_n$ for all $n\ge 0$, we can form
the martingale
$$y_n=x_0+ \sum_1^n \xi_{k-1} (x_k-x_{k-1}).$$

Then,    if   $(x_n)$ is a  martingale which converges
in $L^p(M,\tau)$ ($1<p<\infty$),  if the sequence $\xi=(\xi_n)$ 
is bounded in $\M$  and if $\xi_{n-1}$ commutes with 
$\M_n$ for all $n$,
the transformed martingale $(y_n)$ also converges in $L^p(M,\tau)$. 

\n Indeed, by duality, it suffices to check this for $p\ge 2$,
 and then it is an easy consequence
of (0.1). Note however that the
preceding statement can fail if
one does not assume  that $\xi_{n-1}$ commutes with $\M_n$.
In the case $p\ge 2$, it suffices to assume that 
$\xi_{n-1}$ commutes with $x_n-x_{n-1}$ for all $n$.
 The latter assumption
is used to show that  if, say $\|\xi_{n-1}\|\le 1$, we have 
$(y_n-y_{n-1})(y_n-y_{n-1})^*\le (x_n-x_{n-1})(x_n-x_{n-1})^*$. 
Of course, this assumption can be
relaxed further, all that is needed is to be able to compare 
the ``square functions" associated
to
$(y_n)$ and $(x_n)$ appearing on the right in (0.1).

In section 2 the above inequalities (0.1) and (0.2) are proved. The key point
of our proof
is the following passage:
assuming the above inequalities for some
$1<p<\infty$, then we deduce them for $2p$. The rest of the proof can be
accomplished by iteration (starting from $p=2$),
interpolation and duality. We would like to emphasize that this proof is
entirely  self-contained.

The style of proof of (0.1) and (0.2) is rather old fashioned: it is 
reminiscent of Marcel Riesz's
classical argument for the boundedness of the Hilbert transform on $L^p$
($1<p<\infty$), 
  and
also
of Paley's proof of (0.1)  in  the classical dyadic case ([Pa]), \ie when
$\M_n=L_\infty(\{-1,+1\}^n)$.
 It has been known for many years that  Marcel Riesz's
 argument could be easily adapted to prove the boundedness of the Hilbert
transform
on the vector valued $L^p$-space ($p\ge2$)  $L^p(X)$, when the Banach space
$X$ is 
the Schatten $p$-class $S_p$, or  a non-commutative 
$L^p$-space associated to a trace (the first author learned this from P.
Muhly back in 1976).
 More recently, Bourgain ([B1]) used this to show the unconditionality of
martingale
differences with values
in $S_p$. In other words, he showed that $S_p$
is a UMD space, in the terminology of [Bu2]. (See [BGM] for the case of more
general 
non-commutative 
$L^p$-spaces.)    Recall that a Banach space $X$ is called a UMD space
if, for any $1<q<\infty$, there is a constant $C$ such that, for any
$q$-integrable $X$-valued 
finite martingale
$(x_n)$ on a probability space $(\Omega,{\cal A}, P)$ and for any choice of
sign 
$\epsilon_n=\pm1$, we have (here we write briefly $L_q(X)$ instead of
${L_q(\Omega,{\cal A},
P;X)}$)
$$\|\sum \vp_n (x_n-x_{n-1})\|_{L_q(X)}
\le C\|\sum  
 x_n-x_{n-1} \|_{L_q(X)}=C\sup_n \|x_n\|_{L_q(X)} .\leqno(0.3)$$
We will denote by $C_q(X)$ the best constant $C$ satisfying this.
By well known stopping time arguments (the so-called
``good $\lambda$ inequalities", see [Bu1]) it suffices to have this for
{\it some} $1<q<\infty$, 
for instance for $q=2$ say, and there is a positive constant $K_q$ depending
only on $q$
such that for all $1<q<\infty$
$$K^{-1}_q C_2(X) \le C_q(X)\le K_q C_2(X).\leqno(0.4)$$
Of course, when $X$ is a non-commutative
$L^p$-space, the choice  of $q=p$ gives a nicer form to (0.3).
 The reader is referred to [Bu2] for more information on UMD spaces.

The fact that non-commutative
$L^p$-spaces are UMD ([B1-2, BGM]), which is of course a corollary of our
main result,     can
also be used to prove, by some kind of transference argument, several special
cases of it. This
is explained in section 3. However, although it seems to give better behaved
constants
(when $p\to \infty$),  we do not see how to use this
transference
 idea in the   situation of an arbitrary filtration, as treated in section 2.

In section 3 we give three examples. They are respectively the tensor
products, Clifford algebras and
algebras of free groups. For all of them the preceding inequalities admit a
different proof, that we
outline in the tensor product case. Its main idea  is to transfer a
non-commutative martingale to a
{\it commutative} martingale with values in the corresponding
non-commutative $L^p$-space
$L^p(\M,\tau)$, and to use its unconditionality. This alternate method is, in
fact, our first approach to
non-commutative martingale inequalities, as announced in [PX].

Section 4 is devoted to the Ito-Clifford integral. There we apply our main
inequalities to  give a
characterization for an Ito-Clifford integral to be a
$L^p$-martingale via its integrand. This is the Fermionic analogue of the
square function inequality
for the classical Ito integrals. As a consequence, we extend the
Ito-Clifford integral theory in
$L^2$, developed by Barnett, Streater and Wilde, to $L^p$ for all
$1<p<\infty$.

We include an appendix on the non-commutative analogue of the classical
Fefferman duality between
$H^1$ and $BMO$.

\bigskip
 {\bf Acknowledgement:} We are very grateful to Philippe Biane for several 
fruitful
conversations, and also
to Eric Carlen for kindly providing us with a copy of a preliminary version
of [CK].

\bigskip

\centerline{\bf 1. Preliminaries}

Let $\M$ be a finite von Neumann algebra with a normalized faithful trace
$\tau$. For $1\le p \le \infty$ let $L^p(\M, \tau)$ or simply $L^p(\M)$
denote the associated
non-commutative $L^p$-space. Note that if $p=\infty$,
$L^p(\M)$ is just $\M$ itself with the operator norm; also recall that the
norm in $L^p(\M)$ $(1\le
p<\infty)$ is defined as
$$\|x\|_p = (\tau(|x|^p))^{1/p},\qquad x\in L^p(\M),$$ where $$|x| =
(x^*x)^{1/2}$$ is the usual
absolute value of
$x$.

 Let $a = (a_n)_{n\ge 0}$ be a finite sequence in $L^p(\M)$. Define
$$\|a\|_{L^p(\M; l^2_C)}= \| \big(\sum_{n\ge 0}
|a_n|^2\big)^{1/2}\|_p,\quad \|a\|_{L^p(\M; l^2_R)} =
\|\big( \sum_{n\ge 0} |a^*_n|^2\big)^{1/2}\|_p.\leqno(1.1)$$  This gives two
norms on
the family of all finite
sequences in $L^p(\M)$. To see that, denoting by $B(l^2)$ the algebra of
all bounded operators on
$l^2$ with its usual trace tr, let us consider the von~Neumann algebra
tensor product $\M\otimes
B(l^2)$ with the product trace $\tau \otimes \hbox{tr}$.
$\tau\otimes \hbox{tr}$ is a semifinite faithful trace. The associated
non-commutative $L^p$-space is
denoted by
$L^p(\M \otimes B(l^2))$. Now any finite sequence $a = (a_n)_{n\ge 0}$ in
$L^p(\M)$ can be regarded
as an element in $L^p(\M \otimes B(l^2))$ via the following map $$a
\mapsto T(a) = \left(\matrix{a_0&0&\ldots\cr a_1&0&\ldots\cr
\vdots&\vdots\cr}\right),$$ that is, the
matrix of $T(a)$ has all vanishing entries except those in the first column
which are the $a_n$'s.
Such a matrix is called a column matrix, and the  closure in $L^p(\M
\otimes B(l^2))$ of all column matrices  is called the column subspace of
$L^p(\M \otimes B(l^2))$
(when
$p=\infty$, we take the $w^*$-closure of all column matrices). Then
$$\|a\|_{L^p(\M; l^2_C)} = \|\,
|T(a)|\,
\|_{L^p(\M\otimes B(l^2))} = \|T(a)\|_{L^p(\M \otimes B(l^2))}.$$
Therefore, $\|\cdot\|_{L^p(\M;
l^2_C)}$ defines a norm on the family of all finite sequences of
$L^p(\M)$. The corresponding completion (for $1\le p <\infty$) is a Banach
space, denoted by $L^p(\M;
l^2_C)$. Then $L^p(\M; l^2_C)$ is isometric to the column subspace of
$L^p(\M \otimes B(l^2))$. For
$p=\infty$ we let
$L^\infty(\M; l^2_C)$ be the Banach space of sequences in
$L^\infty(\M)$ isometric by the above map $T$ to the column subspace of
$L^\infty(\M \otimes
B(l^2))$. It is easy to check that a sequence $a = (a_n)_{n\ge 0}$ in
$L^p(\M)$ belongs to $L^p(\M; l^2_C)$ iff
$$\sup_{n\ge 0}\|\big(\sum^n_{k=0} |a_k|^2\big)^{1/2}
\|_p < \infty;$$  if this is the case,
$\big(\sum\limits^\infty_{k=0} |a_k|^2\big)^{1/2}$ belongs to $L^p(\M)$ and
$\big(\sum\limits^n_{k=0}
|a_k|^2\big)^{1/2}$ converges to it in $L^p(\M)$ (relative to the
$w^*$-topology for $p=\infty$).

Similarly (or passing to adjoints), we may show that $\|\cdot\|_{L^p(\M;
l^2_R)}$ is a norm on the
family of all finite sequences in $L^p(\M)$. As above, it defines a Banach
space $L^p(\M; l^2_R)$,
which now is isometric to the row subspace of $L^p(\M
\otimes B(l^2))$ consisting of matrices whose non-zero entries lie only in
the first row.

Observe that the column and row subspaces of $L^p(\M \otimes B(l^2))$ are
1-complemented subspaces.
Therefore, from the classical duality between $L^p(\M \otimes B(l^2))$ and
$L^q(\M
\otimes B(l^2))$ $\big({1\over p} + {1\over q}=1$, $1\le p <\infty\big)$ we
deduce that $$L^p(\M;
l^2_C)^* = L^q(\M; l^2_C)\quad \hbox{and}\quad L^p(\M; l^2_R)^* = L^q(\M;
l^2_R).$$ This
complementation also shows that the families $\{L^p(\M; l^2_C)\}$ and
$\{L^p(\M; l^2_R)\}$ are two interpolation scales, say, for instance,
relative to the complex
interpolation method.

\n Note that, for any finite sequence $ (a_n)_{n\ge 0}$  in $L^p(\M)$, we
have, using tensor
product notation and denoting again by $\|.\|_p$ the norm in ${L^p(\M \otimes
B(l^2))}$
$$\|(\sum a^*_na_n)^{1/2} \|_p=\|\sum a_n\otimes e_{n1}\|_{p}   \quad
\hbox{and}\quad \|(\sum a_na_n^*)^{1/2}\|_p=\|\sum a_n\otimes e_{1n}\|_{p}
.$$

The  following is an extension of a non-commutative version of  H\"older's
inequality from 
[LP], which can be established (perhaps at the cost of an extra factor 2) 
by arguing  as in  
[LP]. For completeness, we include a direct elementary proof 
(without any extra factor) based on the three lines lemma.

\proclaim Lemma 1.1. Let $2\le p\le\infty$. For any finite sequence
$a=(a_n)_{n\ge0}$ in
$L^{2p}(\M)$ and any $A\in L^{2p}(\M)$ we set $B(a,A)=(a_nA)_{n\ge0}$. Then
$$\|B(a,A)\|_{L^p(\M;l^2_R)}\le  \max\left\{ \|a\|_{L^{2p}(\M;l^2_C)}\ , 
\   \|a\|_{L^{2p}(\M;l^2_R)} \right\}\, \|A\|_{2p}\, .\leqno(1.2)$$

\pf By definition, the left side of (1.2) is equal to 
$\|\sum a_n AA^*a_n^*\|^{1/2}_{p/2}$
and, on the other hand, by duality, we have
$$ \|\sum a_n AA^*a_n^*\|_{p/2}=\sup |\psi(B)| \leqno(1.3)$$
with 
$$\psi(B)=\tau(\sum a_n AA^*a_n^* B)$$
and where  the supremum in (1.3) runs over the set of all $B\ge 0$ in $\M$
such that
$\tau(B^r)\le 1$ with $r$ conjugate to $p/2$, or equivalently
with $1/r=1-2/p$.

\n We will apply the three lines lemma to the analytic function $F$ defined
for
$0\le \Re(z) \le 1$ by
$$F(z)=\tau \left(\sum a_n \ (AA^*)^{zp/p'}\ a_n^* \ B^{(1-z)r/p'}\right).$$
Let $\theta=p'/p$ so that $1-\theta =p'/r$. Note that $0\le \theta\le 1 $ and
$F(\theta)=
\psi(B)$. Hence, by the three lines lemma, we have
$$|\psi(B)| = |F(\theta)|\le (\sup_{t\in \reel} |F(it)|)^{1-\theta}
(\sup_{t\in
\reel}|F(1+it)|)^{\theta}.\leqno(1.4)$$
But, by an easy application of H\"older's inequality, we have
$$\sup_{t\in \reel} |F(it)|\le \sup\{\|\sum a_n U a_n^*\|_p \   | \
U\in \M, \ \|U\|\le
1\}, \leqno(1.5)$$ and since $\tau$ is a trace, we also find
$$ \sup_{t\in \reel} |F(1+it)|\le \| (AA^*)^{p/p'}\|_{p'} 
\sup\{\|\sum a_n^* U a_n\|_p \  |
\ U\in \M, \  \|U\|\le 1\}.\leqno(1.6)$$

\n Note that, if $\|U\|\le 1$, we  have 
 $\|\sum a_n U a_n^*\|_p\le \|\sum a_n  
a_n^*\|_p$, and similarly with $a_n^*$ instead of $a_n$. 
Indeed, $\|\sum a_n U a_n^*\|_p= \|(\sum a_n U \otimes e_{1n}) (\sum a^*_n 
\otimes e_{n1})\|_p$,
hence
$$\|\sum a_n U a_n^*\|_p\le \|\sum a_n U \otimes e_{1n}\|_{2p} \| \sum a^*_n 
\otimes
e_{n1}\|_{2p}=\|\sum a_n  a_n^*\|_p.$$
Therefore the inequalities
(1.4), (1.5) and (1.6) combined with (1.3) 
immediately yield the announced result (1.2).
\qed

\n{\bf Remark 1.2.}\ The following example shows
that the right side of (1.2) cannot be simplified
too much: let $\M$ be the algebra of all $N\times N$ complex matrices
equipped with its usual trace, let $A=e_{11}$ and let $a_n=e_{n1}$ for
$n=1,...,N$. Then
$(\sum a_n AA^* a_n^*)^{1/2}=\sum_1^N e_{nn}=(\sum a_n a_n^*)^{1/2}$ and 
$\sum a_n^*a_n = N e_{11}$  so that
$\|(\sum a_n AA^* a_n^*)^{1/2}\|_p=N^{1/p}$, $\| A\|_{2p}=1$  and 
$\|(\sum a_n^*   a_n)^{1/2}\|_{2p}=N^{1/2}$. 
Thus, if $2\le p<\infty$,  for no constant $C$ can the inequality
$ \|B(a,A)\|_{L^p(\M;l^2_R)}\le     
  C\|a\|_{L^{2p}(\M;l^2_R)}  \, \|A\|_{2p}\,  $ be true.
This example also shows that (1.2) fails for $p<2$.
Similarly, the inequality
$ \|B(a,A)\|_{L^p(\M;l^2_R)}\le     
  C\|a\|_{L^{2p}(\M;l^2_C)}  \, \|A\|_{2p}\,  $ 
also fails if $2<p\le \infty$ (take $A=1$ and $a_n=e_{1n}$).

We now  turn to the description of non-commutative martingales and their
square
functions. Let $(\M_n)_{n\ge 0}$
be an increasing sequence of von~Neumann subalgebras of $\M$ such that
$\bigcup\limits_{n\ge 0} \M_n$
generates
$\M$ (in the $w^*$-topology). $(\M_n)_{n\ge 0}$ is called a filtration of
$\M$. The restriction of $\tau$ to $\M_n$ is still denoted by $\tau$. Let
$\E_n = \E(\cdot|\M_n)$ be the conditional expectation of $\M$ with respect
to $\M_n$. $\E_n$ is a
norm 1 projection of $L^p(\M)$ onto $L^p(\M_n)$ for all $1\le p \le
\infty$, and $\E_n(x)\ge 0$
whenever $x\ge 0$. A non-commutative $L^p$-martingale with respect to
$(\M_n)_{n\ge 0}$ is a sequence
$x = (x_n)_{n\ge 0}$ such that $x_n\in L^p(\M_n)$ and
$$\E_m(x_n) = x_m,\qquad \forall \; m=0,1,...,n.$$ Let $\|x\|_p =
\sup\limits_{n\ge 0} \|x_n\|_p$. If
$\|x\|_p<\infty$, $x$ is said to be bounded.
\smallskip \n {\bf Remark 1.3.} Let $x_\infty\in L^p(\M)$. Set
$x_n=\E_n(x_\infty)$ for all $n\ge0$. Then
$x=(x_n)$ is a bounded $L^p$-martingale and $\|x\|_p=\|x_\infty\|_p$;
moreover, $x_n$ converges to
$x_\infty$ in $L^p(\M)$ (relative to the $w^*$-topology in the case
$p=\infty$). Conversely, if
$1<p<\infty$, every bounded $L^p$-martingale converges in $L^p(\M)$, and so
is given by some
$x_\infty\in L^p(\M)$ as previously. Thus one can identify the space of all
bounded
$L^p$-martingales with
$L^p(\M)$ itself in the case $1<p<\infty$.
\smallskip

Let $x$ be a martingale. Its difference sequence, denoted by $dx =
(dx_n)_{n\ge 0}$, is defined as
(with $x_{-1} = 0$ by convention) $$dx_n = x_n-x_{n-1},\qquad n\ge 0.$$ Set
$$S_{C,n}(x) = \left(\sum^n_{k=0} |dx_k|^2\right)^{1/2} \quad \hbox{and}
\quad S_{R,n}(x) = \left(\sum^n_{k=0} |dx^*_k|^2\right)^{1/2}.$$ By the
preceding discussion $dx$
belongs to $L^p(\M; l^2_C)$ (resp.\
$L^p(\M; l^2_R)$) iff $(S_{C,n}(x))_{n\ge 0}$ (resp.\ $(S_{R,n}(x))_{n\ge
0}$) is a bounded sequence
in $L^p(\M)$; in this case,
$$S_C(x) = \left(\sum^\infty_{k=0} |dx_k|^2\right)^{1/2} \quad
\hbox{and}\quad S_R(x) = \left(\sum^\infty_{k=0} |dx^*_k|^2\right)^{1/2}$$
are elements in $L^p(\M)$.
These are the non-commutative analogues of the usual square functions in
the commutative martingale
theory. It should be pointed out that one of $S_C(x)$ and $S_R(x)$ may
exist as element of
$L^p(\M)$ without the other making sense; in other words, the two sequences
$S_{C,n}(x)$ and
$S_{R,n}(x)$ may not be bounded in $L^p(\M)$ at the same time.

Let $1\le p <\infty$. Define $\H^p_C(\M)$ (resp.\ $\H^p_R(\M)$) to be the
space of all
$L^p$-martingales $x$ with respect to $(\M_n)_{n\ge 0}$ such that $dx\in
L^p(\M; l^2_C)$ (resp.\
$dx\in L^p(\M; l^2_R)$), and set
$$\|x\|_{\H^p_C(\M)} = \|dx\|_{L^p(\M; l^2_C)} \quad \hbox{and}\quad
\|x\|_{\H^p_R(\M)} = \|dx\|_{L^p(\M; l^2_R)}.$$ Equipped respectively with
the previous norms,
$\H^p_C(\M)$ and
$\H^p_R(\M)$ are Banach spaces. Note that if $x\in \H^p_C(\M)$,
$$\|x\|_{\H^p_C(\M)} = \sup_{n\ge 0} \|S_{C,n}(x)\|_p = \|S_C(x)\|_p\,,$$
and similar equalities hold
for $\H^p_R(\M)$. Then we define the Hardy spaces of non-commutative
martingales as follows:\ if
$1\le p<2$,
$$\H^p(\M) = \H^p_C(\M) + \H^p_R(\M)$$ equipped with the norm
$$\|x\|_{\H^p(\M)} = \inf\{\|y\|_{\H^p_C(\M)} + \|z\|_{\H^p_R(\M)}\colon \
x=y+z,\quad
y \in \H^p_C(\M), z\in
\H^p_R(\M)\};\leqno(1.7)$$ and if $2\le p <\infty$,
$$\H^p(\M) = \H^p_C(\M)\cap \H^p_R(\M)$$ equipped with the norm
$$\|x\|_{\H^p(\M)} = \max\{\|x\|_{\H^p_C(\M)},\quad
\|x\|_{\H^p_R(\M)}\}.\leqno(1.8)$$ The reason that we have defined
$\H^p(\M)$ differently according to $1\le p<2$ or $2\le p<\infty$ will
become clear in the next
section, where we will show that
$\H^p(\M) = L^p(\M)$ with equivalent norms for all
$1<p<\infty$.

\bigskip

\centerline {\bf 2. The main result}

  In this section $(\M,\tau)$ always denotes a finite von Neumann algebra
equipped with a normalized
faithful trace, and  $(\M_n)_{n\ge0}$ an increasing  filtration of
subalgebras of
$\M$ which generate $\M$.  We keep all notations introduced in the last
section.

In the sequel $\alpha_p$, $\beta_p$, {\it etc}, denote positive constants
depending only on $p$. The
following is the main result of this paper.

\proclaim Theorem 2.1. Let $1<p<\infty$. Let $x = (x_n)_{n\ge 0}$ be an
 $L^p$-martingale with respect to $(\M_n)_{n\ge0}$. Then $x$ is bounded in
$L^p(\M)$ iff $x$ belongs
to $\H^p(\M)$; moreover, if this is the case,
$$\alpha_p^{-1}\|x\|_{\H^p(\M)} \le \|x\|_p \le \beta_p \|x\|_{\H^p(\M)}.
\leqno (BG_p)$$

\smallskip Identifying bounded $L^p$-martingales with their limits, we may
reformulate Theorem~2.1 as
follows.

\proclaim Corollary 2.2. Let $1<p<\infty$. Then $\H^p(\M)=L^p(\M)$ with
equivalent norms.

Corollary 2.2 explains why we have defined, in (1.7) and (1.8), the space
$\H^p(\M)$ and its norm
differently for
$p$ in [1,2) and $[2,\infty)$.
One should note that such a different behavior in the non-commutative case
already appears in the
non-commutative Khintchine inequalities obtained by F.~Lust-Piquard, which
we will recall later on.

Before proceeding to the proof of Theorem 2.1, let us biefly explain our
strategy. Firstly, we prove
the implication ``$(BG_p) \Longrightarrow  (BG_{2p})$'' (this is the key
point of the proof). Then by
iteration (noting that ($BG_2$) is trivial) and interpolation we deduce
($BG_p$) for all $2\le
p<\infty$. Finally,  duality yields ($BG_p$) for $1<p<2$. This is a
well-known approach to the
classical Burkholder-Gundy inequalities in the {\it commutative} martingale
theory. However, in order
to adapt it to the non-commutative setting, one encounters 
several substantial
difficulties.
 Perhaps the main one is the lack of a
reasonable maximal
function in the  non-commutative case. (Note that all the truncation
arguments that appeal to
stopping times appear unavailable or inefficient.)

In the course of the proof we will show (and also need) the following
result, which is the
non-commutative analogue of a classical inequality due to Stein [St].
(See also [B1, Lemma 8] for a similar result in the case of commutative
martingales
with values in a UMD space.) 

\proclaim Theorem 2.3. Let $1<p<\infty$. Define the map $Q$ on all finite
sequences $a = (a_n)_{n\ge
0}$ in $L^p(\M)$ by $Q(a) = (\E_na_n)_{n\ge0}$. Then
$$\|Q(a)\|_{L^p(\M;l^2_C)}\le \gamma_p \|a\|_{L^p(\M;l^2_C)}\,,\quad
\|Q(a)\|_{L^p(\M;l^2_R)}\le
\gamma_p \|a\|_{L^p(\M;l^2_R)}\,.\leqno(S_p)$$ Thus $Q$ extends to a
bounded projection on $L^p(\M;
l^2_C)$ and $L^p(\M; l^2_R)$; consequently, $\H^p(\M)$ is complemented in
$L^p(\M; l^2_C) + L^p(\M,
l^2_R)$ or $L^p(\M; l^2_C) \cap L^p(\M; l^2_R)$  according to $1<p\le 2$ or
$2\le p<\infty$.

 \n {\bf Remark 2.4.}    The inequalities ($BG_p$) imply that  all martingale
difference sequences are unconditional
in
$L^p(\M)$, \ie there is a positive constant $\beta'_p$ such that  for  all
finite martingales
$x$ in 
$L^p(\M)$ we have
$$\|\sum_n\vp_ndx_n\|_p\le
\beta'_p  \|\sum_ndx_n\|_p\,,\quad\forall\;\vp_n=\pm1.\leqno (BG'_p)$$ 
Moreover $\beta'_p \le \alpha_p\beta_p$.

We begin the proof of Theorems 2.1 and 2.3 with some elementary Lemmas.

The inequality  below is well known: indeed, it is a consequence
of the UMD property of
$L^p(\M)$. One can also use the Hilbert transform instead. For the sake of  
completeness, we
will show that it follows from $(BG_p)$. The
following proof is similar to an argument presented in [HP].

\proclaim Lemma 2.5. Let
$\vp=(\vp_n)_{n\ge 0}$ be a sequence of independent random variables on
some probability space
$(\Omega,\F, P)$  such that
$P(\vp_n=1)=P(\vp_n=-1)=1/2$ for all $n\ge 0$. 
Let $\vp'=(\vp'_n)_{n\ge 0}$ be an independent copy of
$\vp$. Let $1<p<\infty$. Suppose  $(BG_p)$. Then for all finite double
sequences
$(a_{ij})_{i,j\ge0}$ in
$L^p(\M)$,
$$\Big(\int_\Omega \|\sum_{0\le i\le j}\vp_i\vp_j'a_{ij}\|_p^p
dP(\vp)dP(\vp')\Big)^{1/p}
\le \alpha_p\beta_p \Big(\int_\Omega \|\sum_{ i,
j\ge0}\vp_i\vp_j'a_{ij}\|_p^pdP(\vp)dP(\vp')\Big)^{1/p}\,.$$

\pf Given $n\ge 0$ let $\F_{2n}$ and $\F_{2n+1}$ be  the
sub-$\sigma$-fields of $\F$
 generated respectively by $\{\vp_0,\cdots,\vp_n\}\cup
\{\vp'_0,\cdots,\vp'_n\}$ and
$\{\vp_0,\cdots,\vp_n,\vp_{n+1}\}\cup \{\vp'_0,\cdots,\vp'_n\}$. Then
$(\F_n)_{n\ge 0}$ is an
increasing filtration of sub-$\sigma$-fields of $\F$. Let  $\exp$ denote
the expectation viewed as a (tracial!) functional   on
$L^\infty(\Omega,\F, P)$. We consider the tensor product
$(\M,\tau)\otimes (L^\infty(\Omega,\F, P),\exp)$ and its increasing
filtration $\M\otimes
L^\infty(\Omega,\F_n, P)$. Hence we have $(BG_p)$ for the corresponding
martingales (noting that
such martingales are in fact {\it commutative} martingales with values in
$L^p(\M)$). Now
given a finite double sequence $(a_{ij})_{i,j\ge0}$ in $L^p(\M)$ we define
a martingale
$f=(f_n)_{n\ge0}$ by
$$f_n={\rm id_\M}\otimes \exp_n\big(\sum_{ i, j\ge0}
\vp_i\vp_j'a_{ij}\big),$$
where $\exp_n$ stands for the conditional expectation of $\F$ with respect
to $\F_n$. Then $(BG_p)$
yields
$$\|\sum_{n\ge 0}\vp''_ndf_n\|_p\le \alpha_p\beta_p\|f\|_p,\quad
\forall\;\vp''_n=\pm1,$$
where the norm $\|\,\cdot\,\|_p$ is understood as it should be, that is, it
is the norm on
$L^p \big(\M\otimes L^\infty(\Omega),\tau\otimes\exp\big)$. Consequently,
$$\|\sum_{n\ge 0}df_{2n}\|_p\le \alpha_p\beta_p\|f\|_p.$$
However,
$$\sum_{n\ge 0}df_{2n}=\sum_{0\le i\le j} \vp_i\vp_j'a_{ij}, $$
whence the announced result.\qed

\proclaim Lemma 2.6. Let $1\le p\le\infty$. Then for all finite sequences
$a=(a_n)_{n\ge0}\subset
L^p(\M)$ we have
$$\|(\sum_{n\ge0}|a_n|^4)^{1/2}\|_p\le \|(\sum_{n\ge0}|a_n|^2)^{1/2}\|_{2p}
\,(\sum_{n\ge0}\|a_n\|^{2p}_{2p})^{1/(2p)}\,.$$

\pf Let $e_{i,j}$  be the matrix in $B(l^2)$ whose entries all vanish but
the one on the position
$(i,j)$ which equals 1. Using the tensor product $\M\otimes B(l^2)$
(already considered in section 1)
we have
$$\eqalign{\|(\sum_{n\ge0}|a_n|^4)^{1/2}\|_p &=\|\sum_{n\ge0}|a_n|^2\otimes
e_{n,0}\|_{L^p(\M\otimes
B(l^2))}\cr &=\|\big(\sum_{n\ge0}a_n^*\otimes
e_{n,n}\big)\big(\sum_{n\ge0}a_n\otimes
e_{n,0}\big)\|_{L^p(\M\otimes B(l^2))}\cr &\le\|\sum_{n\ge0}a_n^*\otimes
e_{n,n}\|_{L^{2p}(\M\otimes
B(l^2))}\|\sum_{n\ge0}a_n\otimes e_{n,0}\|_{L^{2p}(\M\otimes B(l^2))}\cr
&=(\sum_{n\ge0}\|a_n\|^{2p}_{2p})^{1/(2p)}
\|(\sum_{n\ge0}|a_n|^2)^{1/2}\|_{2p}\,.}$$
\qed

In particular, for martingale differences we get the following

\proclaim Lemma 2.7. Let $1\le p\le\infty$. Then for all finite martingales
$x=(x_n)_{n\ge0}\subset
L^{2p}(\M)$ we have
$$\|(\sum_{n\ge0}|dx_n|^4)^{1/2}\|_p\le 2^{1-1/p}
\|x\|_{2p}\|x\|_{\H^{2p}_C(\M)}\,.$$

\pf By Lemma 2.6, it suffices to show
$$\big(\sum_{n\ge0}\|dx_n\|_{2p}^{2p}\big)^{1/(2p)}\le 2^{1-1/p}
\|x\|_{2p}\,.$$ This is trivial for
$p=1$ and $p=\infty$. Then the general case follows by interpolation.\qed

 Now we are prepared to prove Theorems 2.1 and 2.3. The proof  is divided
into several steps.

\n{\bf Proof of Theorems 2.1 and 2.3.}  {\bf Step~1.} ($BG_p$) {\it
implies} ($S_p$).

Let $1<p<\infty$. Suppose ($BG_p$) holds. We will show ($S_p$) holds as
well. 

To this end, fix a
finite sequence $a=(a_k)_{0\le k\le n}\subset L^p(\M)$. We consider 
the tensor product
$(\M, \tau)\otimes (\N, \sigma)$, where $\N=B(l^2_{n+1})$ and
$\sigma=(n+1)^{-1}{\rm tr}$ is the
normalized trace on $B(l^2_{n+1})$. Let $\tilde \E_k=\E_k\otimes {\rm
id}_\N$ denote the conditional
expectation of $\M\otimes \N$ with respect to $\M_k\otimes \N$. Then we have
$(BG_p$) for all martingales relative to the filtration
$(\M_k\otimes\N)_{k\ge 0}$. Now set
$$A_k=(n+1)^{1/p} a_{k}\otimes e_{k,0}\,,\quad 0\le k\le n.$$
Let $\vp=(\vp_n)_{n\ge 0}$ and $\vp'=(\vp'_n)_{n\ge 0}$ be the sequences in
Lemma 2.5. Then, with $\|.\|_p$ denoting here the norm 
in the space ${L^p(\M\otimes \N)}$,  we have
$$\eqalign{  \|Q(a)\|_{L^p(\M;l^2_C)}
&=\|\sum_{k=0}^{n}\tilde\E_k(\vp_kA_k)\|_{p} 
 =\|\sum_{k=0}^{n}\sum_{j=0}^k(\tilde\E_j-\tilde\E_{j-1})(\vp_kA_k)\|_{p}\cr
&=\|\sum_{j=0}^{n}(\tilde \E_j-\tilde
\E_{j-1})\big(\sum_{k=j}^{n}\vp_k A_k\big)\|_{p}\  \cr
\hbox{hence by}\, \,  BG_p \,\, (\cf &  \hbox{Remark  2.4})  \cr &\le
\alpha_p\beta_p
\Big(\int_\Omega\|\sum_{j=0}^{n}\vp'_j(\tilde \E_j-\tilde
\E_{j-1})\big(\sum_{k=j}^{n}\vp_k A_k\big)\|^p_{p}dP(\vp) dP(\vp')\Big)^{1/p}
\cr
\hbox{so by Lemma 2.5,}\,\,\,\,& \cr
& \le (\alpha_p\beta_p)^2 \Big(\int_\Omega\|\sum_{j=0}^{n}\vp'_j(\tilde
\E_j-\tilde
\E_{j-1})\big(\sum_{k=0}^{n}\vp_k A_k\big)\|^p_{p}dP(\vp)
dP(\vp')\Big)^{1/p}\,.  
\cr \hbox{On the other hand}&, \hbox{applying }  (BG_p) \,\, \hbox{once
again, this is}\cr 
 &\le (\alpha_p\beta_p)^3\Big(\int_\Omega\|\sum_{j=0}^{n} (\tilde \E_j-\tilde
\E_{j-1})\big(\sum_{k=0}^{n}\vp_k A_k\big)\|^p_{p}dP(\vp)  \Big)^{1/p}\cr
&=(\alpha_p\beta_p)^3\Big(\int_\Omega\|\sum_{k=0}^{n}\vp_k A_k\|^p_{p}d\vp 
\Big)^{1/p} 
 =(\alpha_p\beta_p)^3 \|a\|_{L^p(\M;l^2_C)}\, .}$$
 Thus, we conclude
$$\|Q(a)\|_{L^p(\M;l^2_C)}\le (\alpha_p\beta_p)^3
\|a\|_{L^p(\M;l^2_C)}\,.$$ Hence $Q$ is bounded on
$L^p(\M;l^2_C)$.  Passing to adjoints yields the boundedness of $Q$  on
$L^p(\M;l^2_R)$. \qed

{\bf Step~2.} {\it  ($BG_p$) implies ($BG_{2p}$).}

Let $1< p<\infty$ and suppose $(BG_p)$. Let $x=(x_n)_{n\ge0}$ be a
martingale in $L^{2p}(\M)$. We
must show $x$ satisfies $(BG_{2p})$. Clearly, we can assume $x$  finite,
that is, there exists
$n\in\nat$ such that $x_k=x_n$ for all $k\ge n$. For simplicity, set
$d_k=dx_k$ (so
$d_k=0$ for all $k\ge n$). Then we write the classical ``Doob identity'':
$$|x_n|^2=x_n^*x_n=S_C(x)^2+\sum_{k\ge0}d_k^*x_{k-1}+
\sum_{k\ge0}x_{k-1}^*d_k\,.\leqno(2.1)$$ 
Hence
$$\eqalign{\|x\|_{2p}^2 &=\||x_n|^2\|_{p}\cr
&\le\|S_C(x)^2\|_{p}+\|\sum_{k\ge0}d_k^*x_{k-1}\|_p+
\|\sum_{k\ge0}x_{k-1}^*d_k\|_p\cr
&=\|x\|^2_{\H^{2p}_C(\M)}+2\|\sum_{k\ge0}d_k^*x_{k-1}\|_p\,.}\leqno(2.2)$$
Observe that
$(d_k^*x_{k-1})_{k\ge0}$ is  a martingale difference sequence. Letting
$y=(y_k)$ be the
corresponding martingale, then by $(BG_p)$, we get
$$\|y\|_p\le \beta_p\|y\|_{\H^p(\M)}\,.\leqno(2.3)$$ Now note that
$$dy_k=d_k^*x_k-d_k^*d_k=\E_k(d_k^*x_n)-|d_k|^2,\quad 0\le k\le
n.\leqno(2.4)$$ Let us first consider
the case $1<p<2$. Then $\|y\|_{\H^p(\M)}\le \|y\|_{\H_C^p(\M)}$, so by
(2.3), (2.4), Lemma 2.7 and
$(S_p)$ (which, by Step~1, holds under $(BG_p)$), we get
$$\eqalign{\|y\|_{\H^p(\M)} &\le\|(\sum_{k=0}^n|d_k|^4)^{1/2}\|_p
         +\|\big(\sum_{k=0}^n|\E_k(d_k^*x_n)|^2\big)^{1/2}\|_p\cr
&\le2^{1-1/p}\|x\|_{2p}\|x\|_{\H_C^{2p}(\M)}
    +\gamma_p\|\big(\sum_{k=0}^nx_n^*d_kd_k^*x_n\big)^{1/2}\|_p\cr &\le
2^{1-1/p}\|x\|_{2p}\|x\|_{\H_C^{2p}(\M)}+\gamma_p\|x\|_{2p}\|x\|_{\H_R^{2p}(
\M)}\cr &\le
(2^{1-1/p}+\gamma_p)\|x\|_{2p}\|x\|_{\H^{2p}(\M)}}\leqno(2.5)$$ If $2\le
p<\infty$, again by (2.3)
and (2.4)
$$\eqalign{\|y\|_{\H^p(\M)} &\le\|(\sum_{k=0}^n|d_k|^4)^{1/2}\|_p\cr
&~~~+\sup\Big\{\|\big(\sum_{k=0}^n|\E_k(d_k^*x_n)|^2\big)^{1/2}\|_p,
\|\big(\sum_{k=0}^n|\E_k(d_k^*x_n)^*|^2\big)^{1/2}\|_p\Big\}.}$$ The first
two terms on the right are
dealt with as before; while by $(S_p)$ and Lemma 1.1, the third term  is
majorized by
$\gamma_p\|x\|_{2p}\|x\|_{\H^{2p}(\M)}$. Thus in the case $2\le p<\infty$,
we have
$$\|y\|_{\H^p(\M)}\le
(2^{1-1/p}+\gamma_p)\|x\|_{2p}\|x\|_{\H^{2p}(\M)}\leqno(2.6)$$ Putting
together (2.2), (2.3), (2.5) and (2.6), we obtain finally
$$\eqalign{\|x\|_{2p}^2 &\le
\|x\|^2_{\H_C^{2p}(\M)}+2\beta_p(2^{1-1/p}+\gamma_p)\|x\|_{2p}\|x\|_{\H^{2p
}(\M)}\cr
&\le\|x\|^2_{\H^{2p}(\M)}+\delta_p\|x\|_{2p}\|x\|_{\H^{2p}(\M)}\,,}$$ where
$\delta_p=2\beta_p(2^{1-1/p}+\gamma_p).$ Therefore, it follows that
$$\|x\|_{2p}\le \beta_{2p}\|x\|_{\H^{2p}(\M)}$$ with $\beta_{2p}={1\over
2}(\delta_p+\sqrt{4+\delta_p^2})$. Thus we have proved the second
inequality of $(BG_{2p}\,)$. The
first one can be obtained in a similar way. Indeed, again by (2.1) and the
previous augument, we get
$$\|x\|_{2p}^2\ge
\|x\|^2_{\H_C^{2p}(\M)}-\delta_p\|x\|_{2p}\|x\|_{\H^{2p}(\M)}\,.$$
Replacing $x_n$
by $x_n^*$  in (2.1), we also have
$$\|x\|_{2p}^2\ge
\|x\|^2_{\H_R^{2p}(\M)}-\delta_p\|x\|_{2p}\|x\|_{\H^{2p}(\M)}\,.$$
Therefore,
$$\|x\|^2_{\H^{2p}(\M)}\le
\|x\|_{2p}^2+\delta_p\|x\|_{2p}\|x\|_{\H^{2p}(\M)}\,,$$ which gives the
first inequality of $(BG_{2p})$.

{\bf Step 3.} {\it $(BG_p)$ for $2\le p<\infty$ and $(S_p)$ for $1<
p<\infty$.}

Evidently, $(BG_2)$ holds with $\alpha_2=\beta_2=1$. Then by Step 2 and
iteration we get
$(BG_{2^n})$ for all positive integers $n$, and so also $(S_{2^n})$ in
virtue of Step 1.

Now we use interpolation to cover all values of $p$ in $[2, \infty)$. This
is easy for $(S_p)$ and
the first inequality of $(BG_p)$. Let us consider, for instance, the first
inequality of $(BG_p)$. By
what we already know  about $(BG_{2^n})$, the linear map $x\mapsto dx$ is
bounded from
$L^{2^n}(\M)$ into $L^{2^n}(\M; l^2_C)$ for every positive integer  $n$. Then
by complex
interpolation, it is bounded from $L^p(\M)$ into $L^p(\M; l^2_C)$ for
$2^n<p<2^{n+1}$, and so for all
$p\in[2,\infty)$. Hence
$$\|dx\|_{L^p(\M; l^2_C)}\le \alpha_p\|x\|_p\,.$$ Passing to adjoints, we
get the same inequality
with $L^p(\M; l^2_R)$ instead of $L^p(\M; l^2_C)$. Thus the first
inequality of $(BG_p)$ holds for
all $2\le p<\infty$. A similar argument applies to
$(S_p)$ for all $2\le p<\infty$. However, the projection $Q$ in Theorem 2.3
is self-adjoint; hence,
we get $(S_p)$ for all $1<p<\infty$, which completes the proof of Theorem
2.3.

Concerning the second inequality of $(BG_p)$, we observe that by duality
and the first inequality of
$(BG_p)$ just proved in $[2, \infty)$, we deduce that  for every $1<p\le2$
and  any martingale $x$ in
$L^p(\M)$ we have
$$\|x\|_p\le \beta_p \inf
\Big\{\|x\|_{\H_C^p(\M)}\,,\|x\|_{\H_R^p(\M)}\Big\}.$$ (Here if $1/p+1/q=1$
(so $2\le q<\infty$),  $\beta_p=\alpha_q$ with $\alpha_q$ being the
constant in the first inequality
of $(BG_q)$; see the next step for more on this). Examining the proof in
Step 2, we see that the
implication ``$(BG_p)
\Longrightarrow (BG_{2p})$'' still holds now  with the help of
$(S_p)$ and the above inequality for all $1<p\le 2$. It  follows that the
second inequality of
$(BG_p)$ holds for all
$2\le p\le 4$. Then Step 2 and iteration   yield the second inequality of
$(BG_p)$ for all $2\le
p<\infty$.

{\bf Step 4.} {\it $(BG_p)$ for $1<p<2$.}

Dualizing $(BG_p)$ in the case $2<p<\infty$, we obtain that if $1<p<2$,
then for all martingales $x$
in $L^p(\M)$
$$\|x\|_p\approx\|dx\|_{L^p(\M;l^2_C)+L^p(\M;l^2_R)}\,.$$ On the other
hand, by Theorem 2.3 (already
proved), $\H^p(\M)$ is complemented in
$L^p(\M;l^2_C)+L^p(\M;l^2_R)$, so the norm of $dx$ in the latter space is
equivalent to the norm of
$x$ in the former.

Therefore, the proof of Theorems 2.1 and 2.3 is now complete.\qed

\n {\bf Remarks.}  {(i)} In Step 3 above, for the proof of the second
inequality
of $(BG_p)$ we have
avoided interpolating the intersection spaces $L^p(\M;l^2_C)\cap
L^p(\M;l^2_R)$ for $p\ge 2$,
although it is shown in [P] that they form an interpolation scale for the
complex method.

 \n { (ii)} The constants $\alpha_p$ and $\beta_p$ given by the above proof
are
not good. In fact, they grow
 exponentially as $p\to\infty$ (see also Remark 3.2 below).

 \smallskip

The inequalities $(BG_p)$ are intimately related to the non-commutative
Khintchine inequalities,
which played an important r\^ole in our first approach to $(BG_p)$ for the
examples considered in the
next section. Let us recall them here for the convenience of the reader.  Let
$\vp=(\vp_n)_{n\ge 0}$ be a sequence of independent random variables on
some probability space
$(\Omega, P)$  such that
$P(\vp_n=1)=P(\vp_n=-1)=1/2$ for all $n\ge 0$.

\proclaim Theorem 2.8 ({\bf Non-commutative Khintchine inequalities,  [LP,
LPP]}).
\hfill\break Let
$1\le p<\infty$. Let $a = (a_n)_{n\ge
0}$ be  a finite sequence  in
$L^p(\M)$.  
 \item{\rm (i)} If $2\le p<\infty$,
$$\|a\|_{L^p(\M;l^2_C)\cap L^p(\M; l^2_R)}\le\Big(\int\limits_\Omega
\|\sum_{n\ge 0}
\vp_na_n\|_p^2dP(\vp)\Big)^{1/2} \le \delta_p\|a\|_{L^p(\M;l^2_C)\cap L^p(\M;
l^2_R)}\,.$$
\item{\rm (ii)} If $1\le p<2$,
$$\alpha\|a\|_{L^p(\M;l^2_C)+ L^p(\M; l^2_R)}\le\Big(\int\limits_\Omega
\|\sum_{n\ge 0}
\vp_na_n\|_p^2dP(\vp)\Big)^{1/2} \le \|a\|_{L^p(\M;l^2_C)+ L^p(\M;
l^2_R)}\,,$$  where  $\alpha>0$ is a
absolute constant.

 This result was first proved in [LP] for
$1<p<\infty$ for the Schatten classes. The general statement
as above (including $p=1$) is contained in  [LPP]. Let us also mention that,
as observed
in [P], a combination of the main result in [LPP] with the
type 2 estimate from [TJ] yields that $\delta_p$ is
of order $\sqrt p$  (the best possible) as $p\to \infty$.  One should
emphasize
that for
$1<p<\infty$  the above non-commutative Khintchine inequalities all
follow from  $(BG_p)$ (with
some worse constants, of course). In that special case however, our proof
essentially reduces to the original one in [LP]. 

\smallskip

\n {\bf Remark 2.9.}  (i)   Note that,  by Theorem 2.8,  the unconditionality
of martingale differences expressed in ($BG'_p$) actually  implies (hence is 
equivalent to)
($BG_p$).  Evidently, ($BG_p$) or ($BG'_p$) is no
longer valid for $p=1$. However, in this case
$p=1$, the second inequality of ($BG_p$) remains true (see the corollary in
the appendix).
Consequently, by the above non-commutative Khintchine inequalities ($p=1$),
we deduce the following
substitute for ($BG'_1$): for any finite martingale $x$ in $L^1(\M)$
$$\sup_{\vp_n=\pm1}\|\sum_n\vp_ndx_n\|_1\approx\|dx\|_{L^1(\M;l^2_C)+L^1(\M;
l^2_R)}\,.$$

(ii) Clearly, $(BG'_p)$ implies the well-known fact (\cf [B1, BGM]) 
that $L^p(\M)$ is a UMD
space for all
$1<p<\infty$ (take $q=p$ in (0.3)).  In particular 
 if $f=(f_n)_{n\ge0}$ is a finite {\it
commutative} martingale
 defined on some probability space with values in $L^p(\M)$, then
 $$\Big(\int
\|\sum_{n\ge0}\vp_n\big(f_n(\omega)-f_{n-1}(\omega)\big)\|_p^pd\omega
\Big)^{1/p}\le
\beta_p'\sup_{n\ge0}\Big(\int\|f_n(\omega)\|_p^pd\omega\Big)^{1/p}\,,\quad
\forall\;\vp_n=\pm1.\leqno(2.7)$$

\bigskip

\centerline {\bf 3. Examples}

In this section, we  give some examples  for which the corresponding
inequalities $(BG_p)$ can be
proved by a different method  from the one given in section 2. The key idea
of this alternate method
is to transfer a non-commutative martingale in $L^p(\M)$ to a {\it
commutative} martingale with
values in $L^p(\M)$. This then enables us to use the unconditionality of
{\it commutative} martingale
differences with values in $L^p(\M)$. (Recall that $L^p(\M)$ is a UMD
space; see Remark 2.9 in
section 2). Although it does not seem suitable in the general case,  this
transference approach might be of
interest in  other situations.  This explains why we will give a sketch of
this second method in the
tensor product case below. Let us also point out that we have first
obtained the non-commutative
martingale inequalities for these examples, before proving the general
Theorem 2.1 (see [PX]).

\smallskip
\noindent{\bf I. Tensor products.}
 Let $(\A_n)$ be a  sequence of hyperfinite  von Neumann  algebras, $\A_n$
being equipped with a
normalized faithful  trace  $\sigma_n$. Let
$$(\M_n,\tau_n)=\bigotimes_{k=0}^n(\A_k,\sigma_k)\quad
\hbox{and}\quad(\M,\tau)=\bigotimes_{k=0}^\infty(\A_k,\sigma_k)$$  be the
tensor products in the
sense of von Neumann algebras. Thus we have an increasing  filtration
$(\M_n)_{n\ge0}$ of subalgebras of
$\M$ which allows us to consider  martingales. Let us reformulate Theorem
2.1 in  this case as
follows.

\proclaim Theorem 3.1. Let $1<p<\infty$ and  $(\M_n)_{n\ge0}$ be as above.
Then $L^p(\M)=\H^p(\M)$
with equivalent norms.

\n {\bf Remark.} A special case of Theorem 3.1 is the one where all $\A_n$'s
are equal to the algebra of
all $2\times 2$ matrices with its normalized trace. Then $\M$ is the
hyperfinite II$_1$ factor, and
$(\M_n)_{n\ge0}$ is its natural filtration.

\n{\bf Sketch of the transference proof of Theorem 3.1.}  It is not hard to
reduce Theorem 3.1 to the
case  where  all
$\A_n$'s are finite dimensional and simple. Thus we will consider this
special case only.  Then let
$\O_n$ be the unitary group  of $\A_n$, equipped with its normalized Haar
measure $\mu_n$ (noting
that since dim$\A_n<\infty$, $\O_n$ is compact). Set
$$(\O,\mu)=\prod_{n\ge0}(\O_n,\mu_n)\,.$$ For $\o=(\o_0,\o_1,\cdots)\in\O$,
we denote by
$\pi_{\o_n}$  the automorphism of $\A_n$ induced by $\o_n$, i.e.
$$\pi_{\o_n}(a)=\o_n^*a\o_n\,,\quad \forall\;a\in\A_n,$$  and we let
$$\pi_\o=\bigotimes_{n\ge0}\pi_{\o_n}\,.$$
 Then $\pi_\o$ is an automorphism of $M$, and extends to an isometry on
$L^p(\M)$ for all $1\le
p\le\infty$.

Now  for $a\in L^p(\M)$ we define
$$f(a,\omega) = \pi_\omega(a),\qquad \forall\; \omega\in\Omega.$$ Then
$f(a,\omega)$ is strongly
measurable as a function from $\Omega$ to $L^p(\M)$ for every $1\le
p<\infty$. Let $\Sigma_n$ be the
$\sigma$-field on
$\Omega$ generated by $(\omega_k)^n_{k=0}$, and   ${\exp}_n =
{\exp}(\,\cdot\, |\Sigma_n)$ the
corresponding conditional expectation. The key point here  is the following
observation:
$${\exp}_kf(a,\omega) = f(\E_k(a),\omega)\quad \hbox{a.e.\ on } \Omega,
\qquad \forall\; k\ge 0,\;\forall\;a\in L^1(\M).$$ (Roughly speaking, the
automorphism $\pi_\omega$
intertwines the two conditional expectations
$\exp_k$ and $\E_k$.) Then let $x$ be  a finite $L^p$-martingale (so  there
is an $n$ such that $x_k
= x_n$ for all $k\ge n $). Let
$f(\omega) = f(x_n,\omega)$ be the function defined above. Then
$({\exp}_kf)_{k\ge 0}$ is a {\it
commutative} martingale on $\Omega$ with values in
$L^p(\M)$, and by the above observation
$$ {\exp}_kf - {\exp}_{k-1}f = \pi_\omega(dx_k),\quad a.e. $$
 Therefore, since   $L^p(\M)$ is a UMD
space (see [B1, B2, BGM]), with constant $C_p=C_p(L^p(\M))$ in (0.3), we have
$$\int\limits_\Omega \|\sum_{k\ge 0}
\vp_k\pi_\omega(dx_k)\|^p_p d\omega \le (C_p)^p
\int\limits_\Omega \|\pi_\omega(x_n)\|^p_p d\omega\,,\quad\forall\;\vp_k =
\pm 1.$$  But
$\pi_\omega$ is an isometry on $L^p(\M)$; hence
$$\|\sum_{k\ge0}\vp_kdx_k\|_p\le C_p\|x\|_p\,,\quad\forall\;\vp_k = \pm
1.$$  Thus we obtain the
unconditionality of martingale differences in $L^p(\M)$, \ie $(BG'_p)$
(defined at
the end of section 2) with $\beta'_p\le C_p$,
which, together with the non-commutative Khintchine inequalities, implies
easily
$(BG_p)$.\qed

\n {\bf Remark 3.2.}\  In this tensor product case (also in the two
following) the
above transference proof
gives better constants $\alpha_p$ and $\beta_p$ in $(BG_p)$ than the
general proof in section 2.
Indeed, by the argument in [B1-2] or [BGM], one can show that the
constant $C_p$ is 
  $O(p^2)$ (resp. $O(1/(p-1)^2)$)  as $p \to\infty$
(resp. $p\to 1$). Note that, when  $  p\ge 2$,   the preceding proof yields
(in the tensor product
case)
  $\alpha_p\le  C_p$ and $\beta_p\le C_p\delta_p$, 
 and when  $1<p\le 2$, $\alpha_p\le \alpha^{-1}C_p\gamma_p$  and $\beta_p\le
C_p$. Actually, a
more careful use of duality yields that for
$p\ge 2$, we still have $\beta_p\le C_p$.
  Therefore,   the preceding  sketch of proof  yields the following estimates
for $\alpha_p$ and
$\beta_p$ in $(BG_p)$: $\alpha_p$ and $\beta_p$ are 
both of order $O(p^2)$ as
$p\to\infty$, and
respectively  of order
    $O((p-1)^{-6})$ and  $O((p-1)^{-2})$ as $p\to 1$. 

\smallskip
\noindent{\bf II. Clifford algebras.} Our second example concerns Clifford
algebras. We take this
opportunity to  give a brief introduction to von Neumann Clifford algebras
and to prepare ourselves
for the next section.  The reader is referred to [PL], [BR], [S] and [C]
for more information on this
subject.

Let $H$ be a complex Hilbert space with a conjugation $J$. Let $\C(H,J)$ or
simply $\C(H)$ denote the
von~Neumann Clifford algebra associated to the
$J$-real subspace of $H$. $\C(H)$ is a finite von Neumann algebra. Let us
briefly describe $\C(H)$
via its  Fock representation.

Denote by
$\Lambda^n(H)$  the
$n$-fold antisymmetric product of $H$, equipped with the canonical scalar
product:
$$\langle u_1\wedge\cdots\wedge u_n, v_1\wedge\cdots\wedge v_n\rangle =\det
(\langle u_k,
v_j\rangle_{1\le k, j\le  n}).$$
$\Lambda^0(H) = {\comp}\1$, where  $\1$ is the vacuum vector.
 The antisymmetric Fock space $\Lambda(H)$ is the direct sum of
$\Lambda^n(H)$: $$\Lambda(H) = \bigoplus_{n\ge 0}
\Lambda^n(H).$$ Given any $v\in H$ the associated creator
$c(v)$ on $\Lambda(H)$ is linearly defined over antisymmetric tensors by
$$c(v) u_1\wedge\cdots\wedge u_n = v\wedge u_1 \wedge \cdots\wedge u_n.$$
$c(v)$ is bounded on $\Lambda(H)$ and $\|c(v)\| = \|v\|$. Its adjoint
$c(v)^*$ is the annihilator
$a(v)$ associated to $v$.  The creators and annihilators satisfy the
following canonical
anticommutation relation (CAR):
 $$\{c(u), a(v)\} = \langle u,v\rangle,\qquad \{c(u), c(v)\}=0,\quad
\forall u, v\in H$$ where
$\{S,T\} = ST+TS$ stands for the anticommutator of $S$ and $T$. The
Fermion field $\Phi$ is then
defined  by
$$\Phi(v) = c(v) + a(Jv),\qquad \forall v\in H.$$
$\Phi$ is a linear map from
$H$ to $B(\Lambda(H))$. Moreover $$\{\Phi(u), \Phi(v)\} = 2\langle u,
Jv\rangle,\qquad \forall u,v\in
H.$$ Therefore, if $u$ and $Jv$ are orthogonal, $\Phi(u)$ and
$\Phi(v)$ anticommute. Notice also that $\Phi(v)$ is hermitian for any
$J$-real vector $v$ (i.e.,
$Jv=v$). Then the von~Neumann Clifford algebra $\C(H)$ is exactly the
 subalgebra of $B(\Lambda(H))$ generated by
$\{\Phi(v)\colon \ v\in H\}$. Observe that if $\{e_i\colon
\ i\in I\}$ is a $J$-real orthonormal basis of $H$,
$\{\Phi(e_i)\colon \ i\in I\}$ is a family of anticommuting hermitian
unitaries, and it generates
$\C(H)$.

The vector state on $B(\Lambda(H))$, given by the vacuum $\1$, induces a
trace $\tau$ on
$\C(H)$:
 $\tau(x) = \langle x(\1),\1\rangle$ for any $x\in \C(H)$. Let $L^p(\C(H))$
denote the associated
non-commutative $L^p$-space.

If $K$ is a $J$-invariant closed subspace  of $H$, $\C(K)$ is naturally
identified as a subalgebra of
$\C(H)$. Now let
$(H_n)_{n\ge 0}$ be an increasing sequence of $J$-invariant closed subspaces
of
$H$ such that $\overline{\bigcup\limits_{n\ge  0}H_n} = H$. Then the
corresponding von~Neumann
Clifford  algebras $(\C(H_n))_{n\ge 0}$ form a filtration of von~Neumann
subalgebras of $\C(H)$. We
will call a non-commutative martingale with respect to $(\C(H_n))_{n\ge 0}$
a Clifford martingale.
Therefore, by Theorem 2.1, we have inequalities $(BG_p)$ for Clifford
martingales. In fact, this
Clifford martingale case can be easily reduced to Theorem 3.1 (the tensor
product case) with the help
of the classical Jordan-Wigner transformation.

\smallskip

Let us consider only a special case for Clifford martingales, where
 $\dim H_n = n$ for all $n\ge 0$. Fix a
$J$-real orthonormal basis $(e_n)_{n\ge 1}$ of $H$ such that $e_n \in H_n
\ominus H_{n-1}$ for all $n\ge 1$. Then $\C_n = \C(H_n)$ is the
$C^*$-algebra generated by $\{\Phi(e_k)\}^n_{k=1}$ and of dimension $2^n$.
For convenience we set
$e_0 =1$ and $e_{-1}=0$. Let $x = (x_n)_{n\ge 0}$ be a Clifford
$L^p$-martingale. Then $dx_n$ can be
written as
$$dx_n = \varphi_n(e_1,\ldots, e_{n-1})\Phi(e_n),$$ where $\varphi_n =
\varphi(e_1,\ldots, e_{n-1})$
belongs to
$L^p(\C_{n-1})$. Let $\varphi = (\varphi_n)_{n\ge 0}$ and $\C=\C(H)$.

\proclaim Proposition 3.3. Let $1\le p\le \infty$ and $x = (x_n)_{n\ge 0}$
be a bounded Clifford
$L^p$-martingale as above. Then $\|dx\|_{L^p(\C; l^2_R)} =
\|\varphi\|_{L^p(\C;l^2_R)}$ and
$${1\over 2}\|\varphi\|_{L^p(\C; l^2_C)} \le \|dx\|_{L^p(\C; l^2_C)} \le
2\|\varphi\|_{L^p(\C;
l^2_C)}.$$

\pf Since $\Phi(e_n)$ is unitary (and hermitian), we have
$\|dx\|_{L^p(\C; l^2_R)} = \|\varphi\|_{L^p(\C; l^2_R)}$.

\n To prove the inequalities on $L^p(\C; l^2_C)$ we need the grading
automorphism (or parity)
$G$ of $\C$: $G$ is uniquely determined by
$$G\big(\Phi(v_1) \ldots\Phi(v_n)\big) = \Phi(-v_1)\ldots\Phi(-v_n),\quad
\forall\;v_k\in H, 0\le k\le n.$$  This means that
$G$ is the automorphism induced by minus the identity of $H$. Recall that
$a\in L^p(\C)$ is called
even (resp.\ odd) if
$G(a) = a$ (resp.\ $G(a) = -a$). We have the decomposition $L^p(\C) =
L^p(\C^+) \oplus L^p(\C^-)$
into even and odd parts; more precisely for any $a\in L^p(\C)$
$$a = {a+G(a)\over 2} + {a-G(a)\over 2} = a^+ + a^-.$$  Since $G$ is
isometric on
$L^p(\C)$,
$$\max(\|a^+\|_p, \|a^-\|_p) \le \|a\|_p \le \|a^+\|_p +
\|a^-\|_p.$$   Now for $x = (x_n)_{n\ge 0}$ as in the proposition we have
$$G(dx_n) =-G(\varphi_n)\Phi(e_n);$$  so $(dx_n)^+ = \varphi^-_n
\Phi(e_n)$. Notice that
$\varphi^-_n\in L^p(\C_{n-1}^-)$. Then by the anticommutation of
$\Phi(e_n)$ with $\Phi(e_k)$ $(1\le k\le n-1)$ we get
$\varphi^-_n\Phi(e_n) = -\Phi(e_n)\varphi^-_n$. Therefore
$(dx_n)^+ = - \Phi(e_n)\varphi^-_n$; hence, since
$\Phi(e_n)$ is unitary,
$$\|(dx_n^+)_{n\ge 0}\|_{L^p(\C;l^2_C)} =
\|(\varphi^-_n)_{n\ge 0}\|_{L^p( \C;l^2_C)}.$$  Similarly,
$$\|(dx^-_n)_{n\ge 0}\|_{L^p(\C; l^2_C)} =
\|(\varphi^+_n)_{n\ge 0}\| _{L^p(\C; l^2_C)}.$$  Combining the preceeding
inequalities, we get
$${1\over 2}
\|\varphi\|_{L^p(\C; L^2_C)} \le \|dx\|_{L^p(\C; l^2_C)}
\le 2\|\varphi\|_{L^p(\C;l^2_C)},$$ proving the proposition.\qed

Let us record explicitly the following  consequence of Theorem~2.1 and
Proposition~3.3.

\proclaim Corollary 3.4. Let $1<p<\infty$ and $x = (x_n)_{n\ge 0}$ be as in
Proposition~3.3. Then if
$2\le p<\infty$ we have
$$\|x\|_{\H^p(\C)} \approx \max\{ \|\varphi\|_{L^p(\C; l^2_C)} ,\
\|\varphi\|_{L^p(\C;l^2_R)}\},$$
and if $1<p<2$ we have
$$\|x\|_{\H^p(\C)} \approx \inf \{\|\varphi'\|_{L^p(\C; l^2_C)} +
\|\varphi''\|_{L^p(\C;l^2_R)}\},$$ where the infimum runs over all
$\varphi'\in L^p(\C; l^2_C)$,
$\varphi''
\in L^p(\C;l^2_R)$ such that $\varphi = \varphi'+\varphi''$ and
$\varphi'_n, \varphi''_n\in L^p(\C_{n-1})$ for all $n\ge 0$.
\smallskip

\noindent{\bf III. Free group  algebras.} Let $\free_n$ be the free group
of $n$ generators.  Let
$vN(\free_\infty)$ be the von Neumann algebra of $\free_n$, equipped with
its standard normalized
trace $\tau$.   $vN(\free_n)$  is naturally identified as a subalgebra of
$vN(\free_{n+1})$, so that $\big(vN(\free_n)\big)_{n\ge1}$ is an increasing
filtration
 of von Neumann subalgebras of
$vN(\free_\infty)$, which generate $vN(\free_\infty)$. For convenience, we
put
$vN(\free_0)=\comp\1$.  Thus we can consider martingales with respect to
$\big(vN(\free_n)\big)_{n\ge0}$. Let
$\H^p\big(vN(\free_\infty)\big)$ denote the corresponding Hardy space.
Then  Theorem 2.1 gives

\proclaim Theorem 3.5.  Let $1<p<\infty$. Then
$$\H^p\big(vN(\free_\infty)\big)=L^p\big(vN(\free_\infty)\big)\quad
\hbox{with}
\;\hbox{equivalent}\;\hbox{norms}.$$

Let us emphasize that, a priori,  the above situation is quite different from
the
one considered in the
tensor product case, since  $vN(\free_n)$ is not hyperfinite as soon as
$n\ge2$. However,  Theorem 3.5
also admits an alternate proof, which appears as a limit case of
  the tensor product
case: indeed, as Philippe  Biane kindly pointed out to us,  this can be done
via random
 matrices with the help of
Voiculescu's limit theorem [V]. We omit
the details. Note that   again this argument yields better 
constants when $p$ tends to infinity,
the same ones as  indicated in Remark 3.2.

\bigskip

\centerline{\bf 4. Applications to the Ito-Clifford integral}

In this section $H$ denotes $L^2({\reel}_+)$ with its usual Lebesgue
measure and complex conjugation;
$\C=\C(H)$ is the associated von~Neumann Clifford algebra equipped with its
normalized trace $\tau$.
For $t\ge 0$ let $H_t$ denote the subspace $L^2(0,t)$ and $\C _t =
\C(H_t)$. Clearly,
$\C_0 = \comp$ and $\C_s\subset \C_t$ for $0\le s\le t$. Let $\E_t =
\E(\cdot \mid \C_t)$ be the
conditional expectation of $\C$ with respect to $\C_t$. Thus we have a
continuous time filtration of
von~Neumann subalgebras
$(\C_t)_{t\ge 0}$ of $\C$, which generate $\C$. All the notions for
discrete martingales in section~1
can be transferred to this continuous time setting. Thus a Clifford
$L^p$-martingale is
 a family $X = (X_t)_{t\ge 0}$ such that $X_t \in L^p(\C_t)$ and
$\E_sX_t=X_s$ for $0\le s\le t$; if additionally $\|X\|_p =
\sup\limits_{t\ge 0} \|X_t\|_p<\infty$,
$X$ is said to be bounded. In this section, unless otherwise stated all
martingales are Clifford
martingales with respect to
$(\C_t)_{t\ge 0}$. The main result here is the analogue of Theorem~2.1 for
these Clifford
martingales. We will deduce it from Theorem~2.1 by discretizing continuous
time Clifford martingales.
This reduction from continuous time to discrete time will be done via the
Ito-Clifford integral
developed by Barnett, Streater and Wilde, who had extended the classical
Ito integral theory to
Clifford
$L^2$-martingales. They showed that any Clifford
$L^2$-martingale admits an Ito-Clifford integral representation. The
Clifford martingale inequalities
below will allow us to extend this Ito-Clifford integral theory from
$L^2$-martingales to
$L^p$-martingales for any
$1<p<\infty$. As a consequence, we will show that any Clifford
$L^p$-martingale $(1<p<2)$ has an
Ito-Clifford integral representation.

Let us first recall the Ito-Clifford integral defined
 in [BSW1-2]. For given
$t\ge 0$ let $\Phi_t = \Phi(\chi_{[0,t)})$ (recalling that $\Phi$ is the
Fermion field defined in
section~3). Then $\Phi_t$ is hermitian and belongs to $\C_t$; by the
canonical anticommutation
relations, $(\Phi_t-\Phi_s)^2 = t-s$ for $0\le s\le t$. $\Phi_t$ is the
Fermion analogue of Brownian
motion.

Like in the classical Ito integral, Barnett, Streater and Wilde develop
their Ito-Clifford integral
by first defining the integrals of simple processes. A simple adapted
$L^p$-process is a function
$f\colon \ {\reel}_+ \to L^p(\C)$  such that $f(t) \in L^p(\C_t)$ for $t\ge
0$ and
$$f(t) = \sum_{k\ge 0} f(t_k) \chi_{[t_k,t_{k+1})}(t),$$ where $(t_k)_{k\ge
0}$ is a subdivision of
${\reel}_+$, i.e., $0=t_0 < t_1 < \cdots$ increasing to $+\infty$. For such
an $f$ we define its
Ito-Clifford integral as follows:\ for $t_k \le t < t_{k+1}$
$$X_t = \int^t_0 f(s) d\Phi_s = \sum^{k-1}_{j=0} f(t_j) (\Phi_{t_{j+1}} -
\Phi_{t_j}) + f(t_k) (\Phi_t-\Phi_{t_k}).$$ Clearly, $X = (X_t)_{t\ge 0}$
is a Clifford
$L^p$-martingale; and if $p=2$,
$$\|X_t\|^2_2 = \int^t_0 \|f(s)\|^2_2ds,\qquad \forall t\ge 0.$$ This
identity allows one to define
the Ito-Clifford integral of any ``adapted $L^2$-process'' $f$ belonging to
$L^2_{\rm loc}({\reel}_+;
L^2(\C))$: $$X_t = \int^t_0 f(s) d\Phi_s,\qquad \forall t\ge 0.$$
$(X_t)_{t\ge 0}$ is again a Clifford $L^2$-martingale and the above
identity still holds. Conversely,
any Clifford $L^2$-martingale admits such an Ito-Clifford integral
representation (cf. [BSW1]).

As in the discrete case, for any simple adapted process $f$ we define
$$S_{C,t}(f) = \big[\int^t_0 f^*(s) f(s)ds\big]^{1/2} \quad \hbox{and}
\quad S_{R,t}(f) = \big[\int^t_0 f(s) f^*(s)ds\big]^{1/2}.$$ Let
$\S^p_{ad}$ be the linear space of
all simple adapted $L^p$-processes and $\S^p_{ad}[0,t]$ its subspace of
processes vanishing in
$(t,\infty)$. Then like in the case of discrete time, $\|S_{C,t}(f)\|_p$ and
$\|S_{R,t}(f)\|_p$ define two norms on $\S^p_{ad}[0,t]$. The completions of
$\S^p_{ad}[0,t]$ with respect to them are denoted respectively by
$\H^p_C[0,t]$ and $\H^p_R[0,t]$ for $1\le p<\infty$. Let us point out that
elements in $\H^p_C[0,t]$
and $\H^p_R[0,t]$ can be regarded as measurable operators in
$L^p(\C_t\otimes B(L^2[0,t]))$ (see section~1 about the column and row
subspaces). Let $\H^p_{C,\rm
loc}({\reel}_+)$  (resp.\ $\H^p_{R,\rm loc}({\reel}_+)$) denote the space
of all functions $f\colon \
{\reel}_+ \to L^p(\C)$ whose restrictions to $[0,t]$ belong to
$\H^p_C[0,t]$ (resp.\ $\H^p_R[0,t])$
for all $t\ge 0$. We call elements in $\H^p_{C,\rm loc}({\reel}_+)$ and
$\H^p_{R,\rm loc}({\reel}_+)$
(measurable) adapted $L^p$-processes. As in the discrete case, we define
$$\H^p[0,t] = \H^p_C[0,t] +
\H^p_R[0,t]\quad \hbox{for}\quad 1\le p<2,$$ and
$$\H^p[0,t] = \H^p_C[0,t] \cap \H^p_R[0,t]\quad
\hbox{for}\quad 2\le p<\infty.$$ We endow $\H^p[0,t]$ with the
corresponding sum or intersection
norm. Similarly, we define $\H^p_{\rm loc}({\reel}_+)$.

Now we can state the main result of this section.

\proclaim Theorem 4.1. Let $1<p<\infty$. Then for any $f\in \H^p_{\rm
loc}({\reel}_+)$ its
Ito-Clifford integral
$$X_t = \int^t_0 f(s)d\Phi_s,\qquad t\ge 0$$ is a well-defined Clifford
$L^p$-martingale and
$$\alpha_p^{-1} \|f\|_{\H^p[0,t]} \le \|X_t\|_p \le \beta_p
\|f\|_{\H^p[0,t]},
\qquad \forall\; t \ge 0.$$

\n {\bf Remarks.} (i) Carlen and Kr\'ee [CK] proved that if $p\le 2$ and if
$f$  is a simple adapted
process,  then  the Ito-Clifford integral $(X_t)$ of $f$ satisfies
$$\|X_t\|_p\le \beta_p\min\Big\{\Big\|\big[\int_0^t
|f(s)|^2\,ds\big]^{1/2}\Big\|_p,\;
\Big\|\big[\int_0^t |f(s)^*|^2\,ds\big]^{1/2}\Big\|_p\Big\}.$$ (This
corresponds essentially to the
second inequality of Theorem 4.1 for $p\le 2$.) From this they  deduced
some sufficient conditions
for the existence of Ito-Clifford integrals. They also proved Theorem 4.1 for
$p=4$ (and mentioned that the same argument works for $p=6$ and 8).

(ii) If $2\le p<\infty$, then
$$\H^p_{\rm loc} ({\reel}_+) \subset L^2_{\rm loc}({\reel}_+; L^2(\C));$$
so adapted $L^p$-processes
are adapted $L^2$-processes. Thus the existence of Ito-Clifford integrals
of adapted
$L^p$-processes $(p\ge 2)$ goes back to [BSW1]. Note also that in the case
$p=2$ the inequalities in Theorem 4.1 become equalities (i.e., $\alpha_2 =
\beta_2=1$). This is the
only case already treated in [BSW1]. If $f$ is an adapted $L^1$-process,
then its Ito-Clifford
integral is also a well-defined Clifford $L^1$-martingale
$X = (X_t)_{t\ge 0}$ and we have
$$\|X_t\|_1 \le \beta_1
\|f\|_{\H^1[0,t]},\qquad \forall t\ge 0 $$ (see the corollary in the
appendix and Remark 2.9). Of
course, the reverse inequality fails this time.

We will reduce Theorem~4.1 to simple adapted processes and then apply
Theorem~2.1. For this reduction
to be successful we have to check two things. The first one is the density
of $\S^p_{ad}[0,t]$ in
$\H^p[0,t]$ (this is trivial for $1\le p\le 2$). The second one is that the
norm of a simple adapted
$L^p$-process $f$ in $\H^p[0,t]$ for $1 <p<2$ is equivalent to
$$\inf\{\|g\|_{\H^p_C[0,t]} + \|h\|_{\H^p_R[0,t]}\colon \ f = g+h,\quad
g,h\in \S^p_{ad}\}.$$ These
will be done by the following lemmas.

\proclaim Lemma 4.2. Let $\sigma = (t_k)^\infty_{k=0}$ be a subdivision of
${\reel}_+$. Define the map $Q_\sigma$ over simple adapted processes by
$$Q_\sigma(f)(t) = {1\over t_{k+1}-t_k} \int^{t_{k+1}}_{t_k} \E_{t_k} f(s)
ds,\qquad t_k \le t <
t_{k+1},\quad t\ge 0.$$ Then for $1<p<\infty$, $Q_\sigma$ extends to a
bounded projection on
$\H^p_C[0,t]$ and $\H^p_R[0,t]$ for all $t\ge 0$.

\pf Suppose $f$ is a simple adapted $L^p$-process:
$$f = \sum_{j\ge 0} f(s_j) \chi_{[s_j,s_{j+1})}.$$ By refining the
subdivision $(s_j)_{j\ge0}$ if
necessary we may assume it is finer than $\sigma$. Then
$$Q_\sigma f = \sum_{k\ge 0} \big[\sum_{j:t_k\le s_j <t_{k+1}} \theta_{k,j}
\E_{t_k}f(s_j)\big]
\chi_{[t_k,t_{k+1})},$$ where
$$\theta_{k,j} = {s_{j+1}-s_j\over t_{k+1}-t_k}\quad \hbox{for}\quad t_k
\le s_j < t_{k+1}.$$ Note that
$$\sum_{j:t_k \le s_j < t_{k+1}} \theta_{k,j} = 1,\qquad
\forall\; k\ge0.$$  Observe also the following elementary and well known
inequality: for any
sequence of operators
$(a_j)$ in $B(H)$ ($H$ being a Hilbert space)  and for any finitely
supported sequence
$(\theta_j)$ with $\theta_j\ge 0$ and $\sum
\theta_j=1$, we have  (in the order of $B(H)$)
$$| \sum \theta_j a_j|^2
\le \sum \theta_j |a_j|^2.$$  (Indeed, for any $h$ in $H$, by convexity of
$\|\,\cdot\,\|^2$, we have
$\|\sum \theta_j a_j h\|^2 \le \sum
\lambda_n \|a_j h\|^2$, whence the desired inequality.) Therefore, for all
$k\ge 0$
$$|\sum_{j:t_k\le s_j <t_{k+1}} \theta_{k,j}\E_{t_k} f(s_j)|^2 \le
\sum_{j:t_k \le s_j<t_{k+1}}
\theta_{k,j} |\E_{t_k}(f(s_j))|^2.$$ Now let $t\ge 0$. Without loss of
generality we assume
$t=t_{n+1}$ for some
$n\ge 0$. Then by Theorem 2.3,
$$\eqalignno{\|Q_\sigma f\|_{\H^p_C[0,t]} &= \Big\|\big[\int^t_0 (Q_\sigma
f(s))^* (Q_\sigma
f(s))ds\big]^{1/2}\Big\|_p\cr &\le
\Big\|\big[\sum^n_{k=0} \sum_{t_k\le s_j < t_{k+1}} (s_{j+1}-s_j)
|\E_{t_k}f(s_j)|^2\big]^{1/2}\Big\|_p\cr
 &\le \beta_p
\Big\|\big[\sum^n_{k=0} \sum_{t_k \le s_j <t_{k+1}} (s_{j+1} -s_j)
|f(s_j)|^2\big]^{1/2}\Big\|_p\cr &=
\beta_p \|f\|_{\H^p_C[0,t]}.}$$ Therefore $Q_\sigma$ extends to a bounded
map (projection) on
$\H^p_C[0,t]$. The same reasoning applies to $\H^p_R[0,t]$.\qed

\proclaim Lemma 4.3. Let $1<p<\infty$ and $f\in \H^p_{C,\rm
loc}({\reel}_+)$ (resp.\ $\H^p_{R,\rm
loc}({\reel}_+)$). Then for all $t\ge 0$
$$\lim_\sigma Q_\sigma f = f \quad \hbox{in}\quad \H^p_C[0,t] \quad
(\hbox{resp.}\ \H^p_R[0,t]),$$
where the limit is taken relative to the subdivision $\sigma = (t_k)_{k\ge
0}$ when
$\sup\limits_{k\ge0} (t_{k+1}-t_k)$ goes to zero.

\pf If $f\in \S^p_{ad}$, then $Q_\sigma f = f$ when $\sigma$ is
sufficiently fine; so the lemma is
true for simple adapted processes. The general case is proved by Lemma~4.2
and the density of
$\S^p_{ad}[0,t]$ in
$\H^p_C[0,t]$ and $\H^p_R[0,t]$.\qed

\proclaim Lemma 4.4. Let $1\le p<\infty$. Then $\S^p_{ad}[0,t]$ is dense in
$\H^p[0,t]$ for all $t\ge
0$.

\pf This is trivial for $1\le p<2$ because $\H^p[0,t] = \H^p_C[0,t] +
\H^p_R[0,t]$ in this case. For $2\le p < \infty$ and $f\in \H^p[0,t]$
Lemma~4.3 implies that
$$\lim_\sigma Q_\sigma f  = f\quad \hbox{in}\quad \H^p[0,t].$$ Thus
$\S^p_{ad}[0,t]$ is also dense in
$\H^p[0,t]$.\qed

\proclaim Lemma 4.5. Let $1<p<\infty$ and $f\in \S^p_{ad}$. Then for all
$t\ge 0$
$$\|f\|_{\H^p_C[0,t] + \H^p_R[0,t]} \approx \inf \{\|g\|_{\H^p_C[0,t]} +
\|h\|_{\H^p_R[0,t]}\},$$ where the infimum is taken over all $g,h\in
\S^p_{ad}[0,t]$ such that $f =
g+h$.

\pf Let $f$ be a simple adapted $L^p$-process defined by a subdivision
$\sigma = (t_k)_{k\ge 0}$:
$$f = \sum_{k\ge 0} f(t_k) \chi_{[t_k, t_{k+1})}.$$ Let $g\in \H^p_C[0,t]$,
$h\in \H^p_R[0,t]$ such
that $f=g+h$ and
$$\|g\|_{\H^p_C[0,t]} + \|h\|_{\H^p_R[0,t]} \le 2\|f\|_{\H^p_C[0,t] +
\H^p_R[0,t]}.$$ Then $f = Q_\sigma f = Q_\sigma g + Q_\sigma h$, and
$Q_\sigma g$,
$Q_\sigma h \in \S^p_{ad}$. By Lemma~4.2
$$\eqalign{\|Q_\sigma g\|_{\H^p_C[0,t]} &\le \beta_p \|g\|_{H^p_C[0,t]},\cr
\|Q_\sigma h\|_{H^p_R[0,t]} &\le \beta_p \|h\|_{\H^p_C[0,t]};}$$ whence the
equivalence in the
lemma.\qed

Now we are ready to show Theorem 4.1.

\noindent{\bf Proof of Theorem 4.1.} First consider the case $2\le
p<\infty$. Let
$f\in \S^p_{ad}$:
$$f = \sum_{k\ge 0} f(t_k) \chi_{[t_k,t_{k+1})}.$$
Then (assuming $t=t_n$)
$$X_t = \sum^{n-1}_{k=0} f(t_k) [\Phi(t_{k+1}) - \Phi(t_k)].$$
Thus $(X_{t_k})^n_{k=0}$ is a finite Clifford $L^p$-martingale with respect
to $(\C(H_{t_k}))^n_{k=0}$. Set
$$d_k = X_{t_{k+1}} - X_{t_k} = f(t_k) [\Phi(t_{k+1}) -
\Phi(t_k)].$$ Then by Theorem 2.1
$$\|X_t\|_p \approx \|\big[\sum^{n-1}_{k=0}
|d_k|^2\big]^{1/2}\|_p + \|\big[\sum^{n-1}_{k=0}
|d^*_k|^2\big]^{1/2}\|_p.$$
Since $\Phi(t_{k+1}) -\Phi(t_k)$ is hermitian and
$$[\Phi(t_{k+1}) - \Phi(t_k)]^2 = t_{k+1} - t_k,$$
we have
$$\eqalign{\sum^{n-1}_{k=0} |d^*_k|^2 &= \sum^{n-1}_{k=0} f(t_k) f(t_k)^*
(t_{k+1}-t_k)\cr
&= \int^t_0 f(s) f(s)^* ds.}$$
On the other hand, since $\chi_{[t_k,t_{k+1})}$ is orthogonal to
$L^2(0,t_k)$ and since $f(t_k)\in \C_{t_k}$, by Proposition~3.3 and its
proof
$$\eqalign{\Big\|\big[\sum^{n-1}_{k=0} |d_k|^2\big]^{1/2}\Big\|_p
&\approx \Big\|\big[\sum^{n-1}_{k=0} f(t_k)^* f(t_k) (t_{k+1}-t_k)\big]
^{1/2}\Big\|_p\cr
&= \Big\|\big[\int^t_0 f(s)^* f(s)ds\big]^{1/2}\Big\|_p\cr
&= \|f\|_{\H^p_C[0,t]}\,.}$$
Therefore, we finally deduce that
$$\|X_t\|_p \approx \max(\|f\|_{\H^p_C[0,t]}, \|f\|_{\H^p_R[0,t]})\,,$$
proving Theorem 4.1 in the case $2\le p < \infty$ for simple adapted
$L^p$-processes. The general adapted $L^p$-processes are treated by
approximation by means of Lemma~4.4.

Now suppose $1<p<2$ and $f\in \S^p_{ad}$. Write
$$f = \sum_{k\ge 0} f(t_k) \chi_{[t_k, t_{k+1})}.$$
Since step functions are dense in $L^2[0,t_k]$, by refining $(t_k)_{k\ge
0}$ if necessary we may assume $f(t_k)$ belongs to the von~Neumann algebra
generated by $\{\Phi(t_{j+1}) - \Phi(t_j)\}^{k-1}_{j=0}$. Let $L_k$ denote
the subspace of $H_{t_k}$ spanned by $\{\chi_{[t_j,t_{j+1})}\}^{k-1}
_{j=0}$. Then $\dim L_k = k$ and $f(t_k) \in \C(L_k)$.
 Let $t=t_n$ for some $n\ge 0$. Then
$$X_t = \sum^{n-1}_{k=0} f(t_k) [\Phi(t_{k+1}) - \Phi(t_k)].$$
Thus $(X_{t_k})^n_{k=1}$ is a finite Clifford martingale relative to
$(\C(L_k))^n_{k=1}$. Applying Corollary~3.4 to $(X_{t_k})^n_{k=1}$ we get
$$\|X_t\|_p \approx \inf \Big\{\|\big[ \sum^{n-1}_{k=0} |a_k|^2
(t_{k+1}-t_k)\big]^{1/2}\|_p + \|\big[ \sum^{n-1}_{k=0}
|b_k^*|^2 (t_{k+1}-t_k)\big]^{1/2}\|_p\Big\}\,,$$
where the infimum runs over all $(a_k)$ and $(b_k)$ such that $a_k+b_k =
f(t_k)$ and $a_k,b_k \in \C(L_k)$ for all $0\le k \le n-1$. Let us show that
the
last infimum is equivalent to $\|f\|_{\H^p[0,t]}$. By Lemma~4.5 (recall
that $f \in \S^p_{ad}$) there are $g,h\in\S^p_{ad}$ such that
$$\|g\|_{\H^p_C[0,t]}  + \|h\|_{\H^p_R[0,t]} \le \beta_p
\|f\|_{\H^p[0,t]};$$
moreover, we may assume that $g$ and $h$ are given by the same subdivision
as $f$. Therefore
$$\|g\|_{\H^p_C[0,t]} = \|\big[ \sum^{n-1}_{k=0} |g(t_k)|^2
(t_{k+1}-t_k)\big]^{1/2}\|_p.$$
Applying Theorem 2.3 to the sequence of conditional expectations
$\{\E\big(\cdot
\mid \C(L_k)\big)\}^n_{k=1}$,  we deduce that
$$\|\big[ \sum^{n-1}_{k=0} |\E\big(g(t_k)\big|
\C(L_k)\big)|^2 (t_{k+1}-t_k) \big]^{1/2}\|_p \le
\beta_p \|g\|_{\H^p_C[0,t]}.$$
The same inequality holds
for $h$ and $\H^p_R[0,t]$ in place of $g$ and
$\H^p_C[0,t]$. Since $f(t_k) \in \C(L_k)$,
$$f =
\sum^{n-1}_{k=0} \Big[\E\big(g(t_k) \big| \C(L_k)\big) + \E\big(h(t_k)\big|
\C(L_k)\big)\Big] \chi_{[t_k,t_{k+1})}.$$
Set $a_k = \E\big(g(t_k)\big|
\C(L_k)\big)$ and $b_k = \E\big(h(t_k)\big| \C(L_k)\big)$ for $0\le
k\le n-1$. Then $f(t_k) = a_k+b_k$ and
$$\eqalign{\|\big[\sum^{n-1}_{k=0} |a_k|^2
(t_{k+1}-t_k)\big]^{1/2} \|_p &\le
\beta_p\|g\|_{\H^p_C[0,t]}\,,\cr
\|\big[\sum^{n-1}_{k=0} |b_k^*|^2
(t_{k+1}-t_k)\big]^{1/2}\|_p &\le
\beta_p\|h\|_{\H^p_R[0,t]}\,.}$$
Thus the desired
equivalence follows, and so $$\|X_t\|_p \approx
\|f\|_{\H^p[0,t]}.$$
Therefore, the inequalities of Theorem 4.1 in the case $1<p<2$
 has been proved for simple adapted
processes. Now let $f\in \H^p_{\rm loc}({\reel}_+)$
$(1<p<2)$. Let $f_n\in \S^p_{ad}[0,t]$ converge to $f$ in
$\H^p[0,t]$. Set $$X^n_t = \int^t_0 f_n(s)ds.$$ Then
$$\|X^n_t-X^m_t\|_p \approx \|f_n-f_m\|_{\H^p[0,t]}.$$
Therefore $X^n_t$ converges to some $X_t$ as $n\to
\infty$. It is clear that $(X_t)_{t\ge 0}$ is a Clifford
$L^p$-martingale and $$\|X_t\|_p \approx
\|f\|_{\H^p[0,t]},\qquad \forall t\ge 0.$$ Also
$(X_t)_{t\ge 0}$ is uniquely determined by $f$. Then we
define the Ito-Clifford integral of $f$ to be $(X_t)_{t\ge
0}$. Hence the proof of Theorem~4.1 is complete.\qed

As a consequence of Theorem~4.1 we get the following Ito-Clifford integral
representation for Clifford $L^p$-martingales $(1<p<\infty)$, which extends
to any $p\in (1,\infty)$ the Barnett-Streater-Wilde representation theorem
for $L^2$-martingales.

\proclaim Theorem 4.6. Let $1<p<\infty$. Then for any Clifford
$L^p$-martingale $(X_t)_{t\ge 0}$ there exists an adapted
$L^p$-process $f\in \H^p_{\rm loc}({\reel}_+)$ such that
$$X_t = X_0 + \int^t_0 f(s)d\Phi_s,\qquad \forall t\ge 0.$$

\pf Let $(X_t)_{t\ge 0}$ be a Clifford $L^p$-martingale.
 Without loss of generality assume $X_0=0$. It suffices to construct the
required adapted process over any interval $[0,T]$. Thus fix $T>0$. For any
subdivision $\sigma$ of $[0,T]\colon\ 0 = t_0 <\cdots < t_n = T$ let
$L_\sigma$ denote the subspace of $H_T =  L^2[0,T]$ spanned by
$\{\chi_{[t_k, t_{k+1})}\}^{n-1}_{k=0}$. Since the union of all $L_\sigma$
is dense in $H_T$, $\C_T = \C(H_T)$ is generated by the union of all
Clifford algebras $\C(L_\sigma)$. It follows that $\bigcup\limits_\sigma
\C(H_\sigma)$ is dense in $L^p(\C_T)$. Therefore there exists a sequence
$(X^n_T)_{n\ge 0}$ of $L^p(\C_T)$ such that $\lim\limits_{n\to \infty}
X^n_T = X_T$ in $L^p(\C_T)$ and such that $X^n_T \in \C(H_{\sigma_n})$ for
some subdivision $\sigma_n$ of $[0,T]$. Let $\sigma_n =
(t^n_k)^{N_n}_{k=0}$. Then $X^n_T$ can be written as
$$X^n_T = \sum^{N_n-1}_{k=0} a_{n,k} [\Phi(t^n_{k+1}) - \Phi(t^n_k)]$$
where $a_{n,k}$ belongs to the $C^*$-algebra generated by
$\{\Phi(t^n_{j+1}) - \Phi(t^n_j)\}^{k-1}_{j=0}$ for all $0 \le k \le N_n$
and $n\ge 0$. Put
$$f_n = \sum^{N_n-1}_{k=0} a_{n,k} \chi_{[t^n_k, t^n_{k+1})}.$$
Then $f_n$ is a simple adapted $L^p$-process and
$$X^n_T = \int^T_0 f_n(s)ds.$$
Therefore, by Theorem 4.1
$$\|X^n_T-X^m_T\|_p \approx \|f_n-f_m\|_{\H^p[0,T]},$$
whence $(f_n)_{n\ge 0}$ is a Cauchy sequence in $\H^p[0,T]$, so it
converges to some adapted $L^p$-process $f\in \H^p[0,T]$. Then clearly
$$X_T = \int^T_0 f(s)ds.$$
This finishes the proof of Theorem 4.6.\qed

\n {\bf Remark.} If we identify a Clifford $L^p$-martingale with the
integrand (adapted $L^p$-process) in its Ito-Clifford integral
representation (this is always possible by Theorem~4.6), then Theorem~4.1
can be reformulated as follows: for any $1<p<\infty$ and any $t\ge0$
$$L^p_0(\C_t) = \H^p[0,t]\quad \hbox{with}\;
\hbox{equivalent}\;\hbox{ norms},$$
where
$$L^p_0(\C_t)=\{X\in L^p(\C_t): \tau(X)=0\}.$$
This equivalence can be extended to the
whole  ${\reel}_+$. Let us say that an adapted $L^p$-process $f$
belongs to $\H^p({\reel}_+)$ if
$$\|f\|_{\H^p({\reel}_+)} = \sup_{t\ge 0} \|f\|_{\H^p[0,t]} <\infty.$$
Then for $1<p<\infty$ a Clifford $L^p$-martingale $X= (X_t)_{t\ge 0}$ is
bounded iff the associated adapted $L^p$-process $f$ belongs to
$\H^p({\reel}_+)$; moreover, in
this case we have
$$\|X\|_p = \sup_{t\ge 0} \|X_t\|_p \approx
\|X_0\|_p+\|f\|_{\H^p({\reel}_+)}.$$
Recall also that $X = (X_t)_{t\ge 0}$ is bounded iff $\lim\limits_{t\to
\infty} X_t = X_\infty$ exists in $L^p(\C)$. Identifying the three objects
$X = (X_t)_{t\ge0}$ with $X_0=0$, $f$ and $X_\infty$, we get that
$\H^p({\reel}_+)=L^p_0(\C)$
with equivalent norms.

\bigskip

\centerline{\bf Appendix.}

In this appendix we consider the non-commutative analogue of the classical
duality between the Hardy
space
$H^1$ and BMO of martingales (see [G]). We will show this duality remains
valid in the
non-commutative case.

Let us go back to the general situation presented in section~1. In all what
follows $(\M,\tau)$
denotes a finite von Neumann algebra with a normalized trace $\tau$, and
$(\M_n)$ an increasing
filtration of von Neumann subalgebras of $\M$, which generate $\M$. Recall
that $\E_n$ denotes the
conditional expectation of $\M$ with respect to $\M_n$. In section~1 we
have introduced the Hardy
spaces
$\H^1_C(\M)$, $\H^1_R(\M)$ and $\H^1(\M)$ of martingales with respect to
$(\M_n)$.

Now let us define the corresponding BMO-spaces.  We set
$$\BMO_C(\M)=\big\{a\in L^2(\M):
\sup_{n\ge0}\|\E_n|a-\E_{n-1}a|^2\|_\infty<\infty\big\},$$ where, as
usual, $\E_{-1}a=0$ (recall $|a|^2=a^*a$). $\BMO_C(\M)$ becomes a Banach
space when equipped with the
norm
$$\|a\|_{\BMO_C(\M)}=\big(\sup_{n\ge0}\|\E_n|a-\E_{n-1}a|^2\|_\infty\big)^{1
/2}.$$ Similarly, we
define $\BMO_R(\M)$, which is  the space of all $a$ such that $a^*\in
\BMO_C(\M)$, equipped with  the
natural norm. Finally, $\BMO(\M)$ is the intersection of these two spaces
$$\BMO(\M)=\BMO_C(\M)\cap \BMO_R(\M)$$ and for any $a\in\BMO(\M)$
$$\|a\|_{\BMO(\M)}=\max\{\|a\|_{\BMO_C(\M)},\|a\|_{\BMO_R(\M)}\}.$$ Notice
that if $a_n=\E_na$, then
$$\E_n|a-\E_{n-1}a|^2=\E_n\big(\sum_{k\ge n}|da_k|^2\big).$$ Note also that
$\E_n|a|^2=\E_{n-1}|a|^2+\E_n|a-\E_{n-1}a|^2$, so that
$\E_n|a-\E_{n-1}a|^2\le
\E_n|a|^2$. Therefore, it follows that
$$\|a\|_{\BMO(\M)}\le \|a\|_\infty\,.\leqno(A_1)$$

For simplicity we will denote $\H^1_C(\M)$, $\BMO_C(\M)$, {\it etc.}
respectively by
$\H^1_C$, $\BMO_C$, etc. We will also adapt the identification between a
martingale and its limit
whenever the latter exists. The result of this appendix is the following
duality.

\proclaim Theorem. We have
$(\H^1_C)^*=\BMO_C$ with  equivalent  norms.  More precisely,
\hfill\break
\indent (i) Every $a\in \BMO_C$ defines a continuous linear functional on
$\H^1_C$ by
$$\f_a(x)=\tau(a^*x),\quad\forall\;x\in L^2(\M).\leqno(A_2)$$
\indent (ii) Conversely, any $\f\in(\H^1_C)^*$ is given as above  by some
$a\in\BMO_C$. Moreover,
$${1\over\sqrt3}\|a\|_{\BMO_C}\le \|\f_a\|_{(\H^1_C)^*}\le
\sqrt2\|a\|_{\BMO_C}.$$
\indent The same duality  holds between $\H^1_R$, $\BMO_R$ and between
$\H^1$, $\BMO$ as well:
$$(\H^1_R)^*=\BMO_R\quad\hbox{and}\quad(\H^1)^*=\BMO.$$

\n {\bf Remark.} In the duality $(A_2)$ we have identified an element $x\in
L^2$ with the martingale
$(\E_nx)_{n\ge0}$. It is evident that this martingale is in $\H^1_C$ and
$$\|x\|_{\H^1_C}\le\|x\|_2.$$ Let us also note that from the discussions in
section 1 the family of
finite martingales is dense in
$\H_C^1$, and so is $L^2$. Of course, the same remark applies to $\H^1_C$
and $\H^1$ as well.

Before proceeding to the proof of the theorem, let us note that the
equivalence constants in (ii)
above are the same as in [G]. In fact, our proof below is modelled on the
one presented in [G],
although one should be careful to some difficulties caused by the
non-commutativity. However, this
time, they are much less substantial than those appearing in the proof of
Theorem 2.1. We will
frequently use the tracial property of
$\tau$ and the following elementary property of expectation:
$$\E_n(abc)=a\E_n(b)c,\quad\forall\;a,c\in\M_n,\;\forall\;b\in\M.$$

\noindent{\bf Proof of the theorem}. (i) Let
$a\in\BMO_C$. Define
$\f_a$ by
$(A_2)$. We must show that
$\f_a$ induces a continuous functional on $\H_C^1$. To that end let $x$ be
a finite $L^2$-
martingale.  Then (recalling our identification between a martingale and
its limit)
$$\f_a(x)=\sum_{n\ge0}\tau(da^*_ndx_n)\,.$$ Set, as in section 1
$$S_{C, n}=\big(\sum_{k=0}^n|dx_k|^2\big)^{1/2}\quad\hbox{and}\quad
S_{C}=\big(\sum_{k=0}^\infty|dx_k|^2\big)^{1/2}.$$ By approximation we may
assume the $S_{C,n}$'s
are  invertible elements in $\M$. Then  by the Cauchy-Schwarz inequality
$$\eqalign{|\f_a(x)|
&=\big|\sum_{n\ge0}\tau(S_{C,n}^{1/2}da^*_ndx_nS_{C,n}^{-1/2})\big|\cr &\le
\big[\tau\big(\sum_{n\ge0}S_{C,n}^{-1/2}|dx_n|^2
S_{C,n}^{-1/2}\big)\big]^{1/2}
\big[\tau\big(\sum_{n\ge0}S_{C,n}^{1/2}|da_n|^2
S_{C,n}^{1/2}\big)\big]^{1/2}\cr
&=\big[\tau(\sum_{n\ge0}S_{C,n}^{-1}|dx_n|^2\big)\big]^{1/2}
\big[\tau\big(\sum_{n\ge0}S_{C,n}|da_n|^2\big)\big]^{1/2}\cr &=I\cdot
II.}$$ We are going to estimate
I and II separately.  First for I we have
$$\eqalign{I^2
&=\sum_{n\ge0}\tau\big([S_{C,n}^2-S_{C,n-1}^2]S_{C,n}^{-1}\big)\cr
&=\sum_{n\ge0}\tau\big([S_{C,n}-S_{C,n-1}][1+S_{C,n-1}S_{C,n}^{-1}]\big)\cr
&\le
\sum_{n\ge0}\tau\big(S_{C,n}-S_{C,n-1}\big)\|1+
S_{C,n-1}S_{C,n}^{-1}\|_\infty\cr &\le
2\tau\big(\sum_{n\ge0}S_{C,n}-S_{C,n-1}\big)\cr &=2\tau
(S_C)=2\|x\|_{\H^1_C},}$$ where we have used
the trivial fact that (noting  $S_{C,n-1}^2\le S_{C,n}^2$)
$$\|S_{C,n-1}S_{C,n}^{-1}\|_\infty^2=\|S_{C,n}^{-1}S_{C,n-1}^2S_{C,n}^{-1}\|
_\infty\le 1.$$ As for
II, set $\theta_0=S_{C,0}$ and  $\theta_n=S_{C,n}-S_{C,n-1}$ for $n\ge1$.
Then
$\theta_n\in\M_n$, and
$$\eqalign{II^2 &=\sum_{n\ge0}\tau\big(S_{C,n}|da_n|^2\big)\cr
&=\sum_{k\ge0}\tau\big[\theta_k\sum_{n\ge k}|da_n|^2\big]\cr
&=\sum_{k\ge0}\tau\big[\theta_k\E_k\big(\sum_{n\ge k}|da_n|^2\big)\big]\cr
&\le\sum_{k\ge0}\tau(\theta_k)\|\E_k\big(\sum_{n\ge
k}|da_n|^2\big)\|_\infty\cr
&\le\|a\|_{\BMO_C}^2\|x\|_{\H_C^1}.}$$ Combining the preceding estimates on
I and II, we obtain, for
any finite $L^2$-martingale $x$
$$|\f_a(x)|\le \sqrt2\|a\|_{\BMO_C}\|x\|_{\H_C^1}.$$ Therefore, $\f_a$
extends to a continuous
functional on $\H_C^1$ of norm$\le\sqrt2\|a\|_{\BMO_C}$.

(ii) Now suppose $\f\in(\H_C^1)^*$. Then by the Hahn-Banach theorem, $\f$
extends to a continuous
functional on
$L^1(\M,l^2_C)$ of the same norm. Thus by the duality (see section 1)
$$(L^1(\M,l^2_C)\big)^*=L^\infty(\M,l^2_C),$$ there exists a sequence
$(b_n)\in L^\infty(\M,l^2_C)$
such that
$$\|\sum_{n\ge0}|b_n|^2\|_\infty=\|\f\|^2\quad \hbox{and }\quad
\f(x)=\sum_{n\ge0}b_n^*dx_n,\quad\forall\;x\in\H_C^1.$$ Let
$a=\sum_{n\ge0}\big(\E_nb_n-\E_{n-1}b_n\big)$ (and so
$da_n=\E_nb_n-\E_{n-1}b_n$). Then $a\in L^2$
and
$$\f(x)=\sum_{n\ge0}da_n^*dx_n=\f_a(x),\quad\forall\;x\in \H_C^1.$$
Therefore, $\f$ is given by
$\f_a$ as in (i). It remains to show $a\in\BMO_C$ and to bound
$\|a\|_{\BMO_C}$ by
$\|\f\|$. This is done as follows. If $k-1\ge n\ge0$,
$$\E_n\big[\E_kb_k^*\E_{k-1}b_k\big]=\E_n\big[\E_{k-1}(\E_kb_k^*\E_{k-1}b_k)
\big]=
\E_n\big[\E_{k-1}b_k^*\E_{k-1}b_k\big];$$ similarly
$$\E_n\big[\E_{k-1}b_k^*\E_{k}b_k\big]=
\E_n\big[\E_{k-1}b_k^*\E_{k-1}b_k\big].$$ It then follows that if $k-1\ge
n\ge0$,
$$\eqalign{\E_n\big[|da_k|^2\big]
&=\E_n\big[(\E_kb_k-\E_{k-1}b_k)^*(\E_kb_k-\E_{k-1}b_k)\big]\cr
&=\E_n\big[\E_kb_k^*\E_kb_k-
\E_{k-1}b_k^*\E_{k-1}b_{k}\big]\cr &\le \E_n\big[\E_kb_k^*\E_kb_k\big]\le
\E_n |b_k|^2.}$$ Hence,
$$\eqalign{\|\E_n|a-\E_{n-1}a|^2\|_\infty &=\|\E_n\sum_{k\ge
n}|da_k|^2\|_\infty\cr
&\le\|\E_n\big[|da_n|^2+\sum_{k\ge n+1}|b_k|^2\big]\|_\infty\cr &\le
3\|\sum_{k\ge0}|b_k|^2\|_\infty\le 3\|\f\|^2;}$$ whence
$$a\in\BMO_C\quad \hbox{and}\quad \|a\|_{\BMO_C}\le \sqrt3\|\f\|.$$ Thus we
have  finished the proof
of the theorem concerning $\H_C^1$ and $\BMO_C$. Passing to adjoints yields
the part on $\H_R^1$ and
$\BMO_R$. Finally, the duality between $\H^1$ and $\BMO$ is obtained by the
classical (and easy) fact
that the dual of a sum is the intersection of the duals.\qed

\proclaim Corollary. Let $x\in \H^1$. Then $x_n$ converges in $L^1$ and
$$\|x\|_1\le \sqrt2\|dx\|_{L^1(\M;l^2_C)+L^1(\M;l^2_R)}\le
\sqrt2\|dx\|_{\H^1}.\leqno (A_3)$$

\noindent {\bf Proof.} Let $x\in \H^1$. By the discussions in section 1,
the finite martingale
$(x_0,\cdots, x_n,x_n,\cdots)$ converges to $x$ in $\H^1$. This, together
with $(A_3)$, implies the
convergence of $x_n$ in $L^1$. Thus it remains to show $(\dagger)$; also it
suffices to show the
first inequality of $(A_3)$ for the second one is trivial. To this end fix
$n\ge0$, and choose
$a\in L^1(\M_n)$ such that $\|a\|_\infty\le1$ and
$\|x_n\|_1=\tau(a^*x_n).$ Put $a_k=\E_k(a)$ for $k\ge0$. Then $a_k=a$ for all
$k\ge n$, and
$$\eqalign{\|x_n\|_1 &=\tau\sum_{k=0}^n da_k^*dx_k =\tau\sum_{k=0}^\infty
da_k^*dx_k\cr
&\le\|dx\|_{L^1(\M;l^2_C)+L^1(\M;l^2_R)}\,\|da\|_{{{L^\infty(\M;l^2_C)\cap
L^\infty(\M;l^2_R)}
\over{(\H^1)^\bot}}}\,.} $$ However, by the preceding theorem
$${{L^\infty(\M;l^2_C)\cap L^\infty(\M;l^2_R)}
\over{(\H^1)^\bot}}=\big(\H^1\big)^*\cong\BMO\,.$$
Therefore, by $(A_1)$
$$\|da\|_{{{L^\infty(\M;l^2_C)\cap L^\infty(\M;l^2_R)}
\over{(\H^1)^\bot}}}\le\sqrt2\|a\|_{\BMO}\le \sqrt2\|a\|_\infty\le
\sqrt2.$$ Combining the previous
inequalities we obtain $(A_3)$, and thus complete the proof of the
corollary.     \qed

\bigskip
\vfill\break

\centerline{\bf References}

\item{[AvW]} L. Accardi, M. von Waldenfels (Eds.) Quantum probability and
applications. (Proc. 1988)
Springer Lecture Notes 1442 (1990).

\item{[BSW1]} C.\ Barnett, R.F.\ Streater, I.F.\ Wilde, The It\^o-Clifford
Integral, J.\ Funct.\
Analysis 48 (1982), 172--212;
\item{} The It\^o-Clifford Integral II - Stochastic differential equation,
J.\ London Math.\ Soc. 27
(1983), 373--384;
\item{} The It\^o-Clifford Integral III -- Markov property of solutions to
Stochastic differential
equation, Commun.\ Math.\ Phys. 89 (1983), 13--17; \item{} The
It\^o-Clifford Integral IV -- A
Radon-Nikodym theorem and bracket processes, J.\ Operator Theory 11 (1984),
255-211.

\item{[BSW2]} C.\ Barnett, R.F.\ Streater, J.F.\ Wilde, Stochastic
integrals in an arbitrary
probability gauge space, Math.\ Proc.\ Camb.\ Phil.\ Soc. 94 (1983),
541--551.

\item{[BGM]} E.\ Berkson, T.A.\ Gillespie, P.S.\ Muhly, Abstract spectral
decompositions guaranteed
by the Hilbert transform, Proc.\ London Math.\ Soc. 53 (1986), 489--517.

\item{[B1]} J. Bourgain, Vector valued singular integrals and the $H^1$-BMO
duality, Probability
theory and harmonic analysis (Chao-Woyczynski, ed.), pp.~1--19; Dekker, New
York, 1986.

\item{[B2]} J. Bourgain, Some remarks on Banach spaces in which martingale
differences are
unconditional, Arkiv f\"or Math. 21 (1983) 163-168.

\item{[BR]} O.\ Bratteli, D.W.\ Robinson, Operator algebras and quantum
statistical mechanics II,
Springer-Verlag, 1981.

\item{[Bu1]} D.\ Burkholder, Distribution function inequalities for
martingales, Ann.
Probability 1 (1973) 19-42.

\item{[Bu2]} D.\ Burkholder, A geometrical characterization of Banach spaces
in which martingale
difference sequences are
 unconditional, Ann.\ Probab. 9 (1981), 997-1011.

\item{[C]} J.M.\ Cook, The mathematic of second quantization, Trans.\
Amer.\ Math.\ Soc. 74 (1953),
222--245.

\item{[CK]} E.A.\ Carlen, P.\ Kr\'ee, On martingale inequalities in
noncommutative stochastic
analysis, Preprint, 1996.

\item{[CL]} E. Carlen and E. Lieb, Optimal hypercontractivity for Fermi
fields and related
non-commutative integration inequalities, Comm. Math. Phys., 155 (1993),
27-46.

\item{[G]} A.M. Garsia, Martingale inequalities, Seminar Notes on Recent
Progress, Benjamin Inc. 1973.

\item{[Gr1]} L. Gross, Existence  and uniqueness of physical ground states,
J. Funct. Analysis, 10
(1972), 52-109.

\item{[Gr2]} L. Gross, Hypercontractivity and logarithmic Sobolev
inequalities for the
Clifford-Dirichlet form,  Duke Math. J. 42 (1975), 383-396.

\item{[HP]}  U. Haagerup and G. Pisier.  Factorization of
analytic functions with values in non-commutative
$L_1$-spaces. 
 Canadian Journal of Math. 41 (1989) 882-906.

\item{[LP]} F.\ Lust-Piquard, In\'egalit\'es de Khintchine dans $C_p$
$(1<p<\infty)$, C.R.\ Acad.\ Sci. Paris 303 (1986), 289--292.

\item{[LPP]} F.\ Lust-Piquard, G.\ Pisier, Noncommutative Khintchine and
Paley inequalities, Arkiv
f\"or Mat. 29 (1991), 241--260.

\item{[M]} P.A. Meyer, Quantum probability  for probabilists.  Lect.\ Notes
Math. 1538,  Springer
Verlag, 1995.

\item{[Mi]}  M. Mitrea, Clifford wavelets, singular integrals and Hardy
spaces.  Lect.\ Notes Math.
1575, Springer Verlag, 1994.

\item{[Pa]} R. E. A. C. Paley, A remarkable series of orthogonal functions
(I), Proc. London Math. Soc. 34 (1932) 241-264.

\item{[P]} G. Pisier,     Non-commutative
 vector valued $L_p$-spaces and completely
$p$-summing maps, (Revised February 97) To appear.

\item{[PX]} G. Pisier, Q. Xu, In\'egalit\'es de martingales non
commutatives, C. R. Acad. Sci. Paris,
323 (1996), 817-822.

 \item{[PR]} R.J.\ Plymen, P.L.\ Robinson, Spinors in Hilbert space.
Cambridge University Press, 1994.

\item{[S]} I.E.\ Segal, Tensor algebra over Hilbert spaces II, Ann.\ Math.
63 (1956), 160--175.

\item{[St]} E. M. Stein, Topics in harmonic analysis related to the
Littlewood-Paley theory,
Princeton University Press, Princeton, N.J., 1970.

\item{[TJ]} N. Tomczak-Jaegermann. The moduli of 
convexity and smoothness and the Rademacher averages of trace class
$S_p$. Studia Math. 50 (1974) 163-182.

\item{[V]} D.V. Voiculescu, Limit laws from random matrices and free
products, Invent.Math., 104
(1991), 201-220.

\item{[VDN]} D.V. Voiculescu, K.J. Dykema, A. Nica, Free Random variables,
CRM Monograph Series,
Vol.1, Centre de Recherches Math\'ematiques, Universit\'e Montr\'eal,

 \end